\def\E{\end{document}}
\documentclass[11pt]{article}
\usepackage{amssymb,amsmath}

\usepackage{CJK}

\begin{document}
\title{}
\date{}
\author{}
\title{\bf Quenching Time Optimal Control for Some Ordinary Differential Equations \footnote{This work was partially supported by  the National Natural
Science Foundation of China under grant 11071036 and the National
Basis Research Program of China (973 Program) under grant
2011CB808002.} }
\date{}
\author{Ping Lin\footnote{  School of
Mathematics \& Statistics, Northeast Normal University, Changchun
130024,  P. R. China. E-mail address:
linp258@nenu.edu.cn.}
\\
\\
 }
\maketitle
\newtheorem{theorem}{Theorem}[section]
\newtheorem{definition}{Definition}[section]
\newtheorem{lemma}{Lemma}[section]
\newtheorem{proposition}{Proposition}[section]
\newtheorem{corollary}{Corollary}[section]
\newtheorem{remark}{Remark}[section]
\renewcommand{\theequation}{\thesection.\arabic{equation}}
\catcode`@=11 \@addtoreset{equation}{section} \catcode`@=12

{\bf Abstract.} This paper concerns some time optimal control problems
of three
different ordinary differential equations in $\mathbb{R}^2$. Corresponding to certain initial data and controls, the
solutions of the systems quench
at finite time. The goal to control the systems is to
minimize the quenching time.  To our best knowledge, the study on  quenching time optimal
control problems has not been touched upon. The purpose of this study is to
obtain the existence and the Pontryagin maximum principles of
optimal controls.  We hope that our methods could hint people to study the same problems  with more general vector fields in $\mathbb{R}^d$ with $d\in \mathbb{N}$. We also wish that our results could be extended to the same issue for parabolic equations.

{\bf Key Words.} quenching time, optimal control, nonlinear ordinary
differential equations

{\bf AMS subject classifications.} 49J15, 34A34

\section{Introduction}
$\; \;\;\;$ In this paper, we study some quenching time optimal control problems of three
different ordinary differential equations in $\mathbb{R}^2$. First of all, some notations will be introduced. We use $\|\cdot\|$ and $<\cdot,\cdot>$ to stand for
the
Euclidean norm  and the inner product of $\mathbb{R}^2$.
For each matrix $D$, we use $D^T$ and $\|D\|$ to denote its
transposition and the operator norm, respectively.
Let $
B(\cdot)\triangleq\begin{pmatrix} b_{11}(\cdot)& b_{12}(\cdot)
\\b_{21}(\cdot)& b_{22}(\cdot)
\end{pmatrix}
$ be a nontrivial matrix-value function in the space $L^{\infty}(0,+\infty;\mathbb{R}^{2\times 2})$.
Write $R^+$ for $[0,+\infty)$.
For each
$\rho_0$  given, we set
$$\mathcal{U}_{ad}=\{u:R^+\rightarrow\mathbb{R}^{2} ;\
u\ \mbox{is Lebesgue measurable,}\ \|u(t)\|\leq\rho_0\ \mbox{for
a.e.}\ t\in R^+\}.$$
Each $u(\cdot)\in \mathcal{U}_{ad}$ can be expressed as  $
u(\cdot)=(u_1(\cdot), u_2(\cdot))^T$. Let
$$b_1(\cdot, u(\cdot))\triangleq b_{11}(\cdot)u_1(\cdot)+b_{12}(\cdot)u_2(\cdot),\;\;\mbox{when}\;\;u(\cdot)\in \mathcal{U}_{ad}
$$
and
$$b_2(\cdot, u(\cdot))\triangleq b_{21}(\cdot)u_1(\cdot)+b_{22}(\cdot)u_2(\cdot),\;\;\mbox{when}\;\;u(\cdot)\in \mathcal{U}_{ad}.
$$
For each $C^1$-function $g: E\subset\mathbb{R}^{2}\rightarrow
\mathbb{R}^{2}$, its derivative will be written as
$$
\begin{pmatrix} \displaystyle\frac{\partial g_1}{\partial y_1} &
\displaystyle\frac{\partial g_2}{\partial y_1} \\
\displaystyle\frac{\partial g_1}{\partial y_2} &
\displaystyle\frac{\partial g_2}{\partial y_2}
\end{pmatrix},
$$
where $g=(g_1,g_2)^T$ with $g_i: E\subset\mathbb{R}^{2}\rightarrow
\mathbb{R}^{1}$, i=1, 2.

The controlled systems under consideration are as follows:
\begin{align}
\label{1.1} \left\{\begin{array}{ll} \displaystyle\frac{dy(t)}{dt}=f(y(t))+B(t)u(t),\ \ \ \   t>0,\\
\displaystyle y(0)=y^0.\end{array}\right.
\end{align}
Here,  $y^0=(y_1^0, y_2^0)^T\in \mathbb{R}^{2}$,   $u(\cdot)\in \mathcal{U}_{ad}$ and
 $f\in \Lambda\triangleq\{f^{(1)}, f^{(2)}, f^{(3)}\}$, where
\begin{align}
(i)\ &f^{(1)}=(f_1^{(1)}, f_2^{(1)})^T:\mathbb{R}^2\setminus\{y=(y_1, y_2)^T\in \mathbb{R}^{2};\ y_1=1\}
\rightarrow\mathbb{R}^2, \nonumber \\&f_1^{(1)}(y)=\frac{y_2}{1-y_1},f_2^{(1)}(y)= y_1+y_2;
\nonumber\\(ii)\ &f^{(2)}=(f_1^{(2)}, f_2^{(2)})^T:\mathbb{R}^2\setminus\{y=(y_1, y_2)^T\in \mathbb{R}^{2};\
y_1^2+y_2^2=1\}\rightarrow\mathbb{R}^2,\nonumber \\&f_1^{(2)}(y)=\frac{y_1}{1-\sqrt{y_1^2+y_2^2}},
f_2^{(2)}(y)=\frac{y_2}{1-\sqrt{y_1^2+y_2^2}};\nonumber\\
(iii)\ &f^{(3)}=(f_1^{(3)}, f_2^{(3)})^T:\mathbb{R}^2\setminus\{y=(y_1, y_2)^T\in \mathbb{R}^{2};\ y_1=1\
\mbox{or}\ y_2=1\}\rightarrow\mathbb{R}^2, \nonumber \\&
f_1^{(3)}(y)=\frac{1}{1-y_2},f_2^{(3)}(y)=\frac{1}{1-y_1}.\nonumber
\end{align}
Let
\begin{align}
\label{1.1-1}K_0=\mbox{esssup}_{s\in R^+}\|B(s)\|  \rho_0.
\end{align}
Define
\begin{align}
S^{f^{(1)}}=&\big\{(z_1,z_2)\in \mathbb{R}^2;\ z_1\in (1-\frac{1}{2K_0},1),\ z_2\in(K_0+\frac{1}{K_0}-1,+\infty)\big\}
\nonumber\\ &\cup \big\{(z_1,z_2)\in \mathbb{R}^2;\ z_1\in (1, 1+\frac{1}{2K_0}),\ z_2\in(K_0+1,+\infty)\big\};
\nonumber\\
S^{f^{(2)}}=&\big\{(z_1,z_2)\in \mathbb{R}^2;\ \|z\|\in (1-\frac{1}{2K_0+1},1)\big\}
\cup \big\{(z_1,z_2)\in \mathbb{R}^2;\ \|z\|\in (1, 1+\frac{1}{2K_0})\big\};
\nonumber\\
S^{f^{(3)}}=&\big\{(z_1,z_2)\in \mathbb{R}^2;\ z_1\in (1-\frac{e^{-3/2}}{2K_0},1),\ z_2\in(1-\frac{e^{-3/2}}{2K_0},1)\big\}
\nonumber\\ &\cup \big\{(z_1,z_2)\in \mathbb{R}^2;\ z_1\in (1, 1+\frac{e^{-3/2}}{2K_0}),\ z_2\in(1, 1+\frac{e^{-3/2}}{2K_0})\big\}.
\nonumber
\end{align}
Since for each $f\in \Lambda$, $f$ is continuously
differential over the domain $S^f$,   it is clear
that given $f\in \Lambda$, $y^0\in S^{f}$  and $u\in \mathcal{U}_{ad}$, the controlled system $(1.1)$
has  a unique solution. We denote this solution  by $y(\cdot\ ;f,y^0,u)
=(y_1(\cdot\ ;f,y^0,u),\ y_2(\cdot\ ;f,y^0,u))^T$, and write $[0,T_{max}{(f,y^0,u)})$ for its maximal interval of
existence.

It is shown in Section 2 that given $f\in \Lambda$, $y^0\in S^f$ and $u\in \mathcal{U}_{ad}$, there exists a
time $T_q(f,y^0,u)$ with $T_q(f,y^0,u)\leq T_{max}{(f,y^0,u)}$ holding the property that
\begin{equation}\label{wang1}
0<T_q(f,y^0,u)<+\infty,\ \lim\limits_{t\rightarrow T_q(f,y^0,u)}\|f(y(t;f,y^0,u))\|=+\infty
\end{equation}
and
\begin{align}
\|f(y(t;f,y^0,u))\|<+\infty\;\;\mbox{as}\;\; t\in [0,T_q(f,y^0,u)).\nonumber\end{align}
We say that the solution $y(\cdot\ ;f,y^0,u)$ quenches at the finite time $T_q(f,y^0,u)$ and $T_q(f,y^0,u)$ is the quenching time of the solution $y(\cdot\ ;f,y^0,u)$.

The purpose of this paper is to study the existence and  the Pontryagin maximum principles for the following time optimal control
problems:
$$({P})_{y^0}^f\ \ \ \min\limits_{\;u\in
\mathcal{U}_{ad}}\big \{ T_{q}(f,y^0, u)\big\},\;\;\mbox{where}\;\; f\in \Lambda\ \mbox{and}\ y^0\in S^f.$$
Because of (\ref{wang1}), for each  $f\in \Lambda$ and $y^0\in S^f$, there exists a number $t^*(f, y^0)$ in $R^+$ such that
$$
t^*(f, y^0)=\inf\limits_{u\in \mathcal{U}_{ad}} T_q(f,y^0,
u),
$$
which is called  the optimal time for the problem
$(P)^f_{y^0}$. A control $u^*(\cdot)$ in the set $\mathcal{U}_{ad}$ holding
the property: $T_q(f,y^0, u^*)=t^*(f,y^0)$, is called an optimal control, while
the solution $y(\cdot\;; f, y^0,u^* )$ is called the optimal state
corresponding to $u^*$ for the problem $(P)_{y^0}^{f}$. We shall simply write
$y^*(\cdot)$ for the optimal state $y(\cdot\ ; f, y^0, u^*)$.

The main results of this paper are as follows:

\begin{theorem} Given  $i\in\{1,2,3\}$ and $y^0\in S^{f^{(i)}}$, the problem $(P)^{f^{(i)}}_{y^0}$ has optimal controls.
\end{theorem}

\begin{theorem}

Let $y^0\in S^{f^{(1)}}$. Then, Pontryagin's maximum principle holds for the problem $(P)^{f^{(1)}}_{y^0}$. Namely,
if $t^*$ is the optimal time, $ u^*$ is an optimal control and $ y^*$ is the corresponding optimal state
for the problem $(P)^{f^{(1)}}_{y^0}$, then there is a nontrivial function $\psi(\cdot)$ in the space
$C([ 0,t^* ]; \mathbb{R}^{2})$ satisfying
\begin{align}\label{1.3}\left\{\begin{array}{ll}
\psi(t)=\displaystyle\int_t^{t^*} f_y^{(1)}(y^*(\tau))\psi(\tau)d\tau \ \
\mbox{for all}\; \ t\in [0,t^* ),\\
\psi(t^*)=0\end{array}\right.
\end{align}
and
\begin{align}\label{1.4}
\max_{  \|u\|\leq\rho_0   }<\psi(t), B(t)u> =<\psi(t),
B(t)u^*(t)>\;\;\mbox{for a.e.} \; t\in [0, t^* ].
\end{align}
Besides, it holds that
\begin{align}\label{1.4-1}
t^*\leq (y^0_1-1)^2.
\end{align}
\end{theorem}

\begin{theorem}
Let $y^0\in S^{f^{(2)}}$. Then, Pontryagin's maximum principle holds for the problem $(P)^{f^{(2)}}_{y^0}$. Namely,
if $t^*$ is the optimal time, $ u^*$ is an optimal control and $
y^*$ is the corresponding optimal state for the problem $(P)^{f^{(2)}}_{y^0}$, then
there is a nontrivial function ${\psi}(\cdot)$ in the
space $C([ 0,t^* ]; \mathbb{R}^{2})$ satisfying
\begin{align}\label{1.5}\left\{\begin{array}{ll}
{\psi}(t)=\displaystyle\int_t^{t^*} f_y^{(2)}(y^*(\tau)){\psi}(\tau)d\tau\ \
\mbox{for all}\; \ t\in [0,t^* ),\\
{\psi}(t^*)=0\end{array}\right.
\end{align}
and
\begin{align}\label{1.6}
\max_{  \|u\|\leq\rho_0   }<{\psi}(t), B(t)u> =<{\psi}(t),
B(t)u^*(t)>\;\;\mbox{for a.e.} \; t\in [0, t^* ].
\end{align}
Besides, it holds that
\begin{align}\label{1.6-1}
t^*\leq \frac{(2K_0+1)(\|y^0\|-1)^2}{2K_0}.
\end{align}
\end{theorem}

It is worth mentioning that Problem $(P)_{y^0}^{f^{(3)}}$ is more complicated than
the other two problems since the target set of Problem $(P)_{y^0}^{f^{(3)}}$ is more complicated than the target sets
of the other two problems. Indeed, the target set of Problem $(P)_{y^0}^{f^{(3)}}$ is
$\{y\in \mathbb{R}^{2};\ y_1=1\ \mbox{or}\ y_2=1\}$; while the target sets for Problem $(P)_{y^0}^{f^{(1)}}$ and Problem $(P)_{y^0}^{f^{(2)}}$ are
accordingly
$\{y\in \mathbb{R}^{2};\ y_1=1\}$ and $\{y\in \mathbb{R}^{2};\ \|y\|=1\}$. This is why we only obtain the existence but not Pontryagin's maximum principle for
Problem $(P)_{y^0}^{f^{(3)}}$.

The concept of quenching was first introduced by H. Kawarada in \cite{Kawarada} for a
nonlinear parabolic equation.  Then there have been  many literatures concerning  the
properties of parabolic  differential equations with
 quenching behavior  (see
\cite{Chan}, \cite{Chan 1}, \cite{B. Hu} and references therein). To our best knowledge, the study on  quenching time optimal
control problems has not been touched upon. In this paper, we focus on  quenching time optimal control problems for ordinary differential equations
with  three particular vector fields $f^{(i)}$, $i=1,2,3,$ in $\mathbb{R}^2$. We hope that our methods could hint people to study the same problems  with more general vector fields in $\mathbb{R}^d$ with $d\in \mathbb{N}$. We also wish that our results could be extended to the same issue for parabolic equations.

The differential systems whose solutions have the behavior of
quenching arise in the study of the electric current transients phenomena in polarized ionic conductors.
It is also significant in the theory of ecology and environmental studies. In
certain cases, the quenching of a solution is desirable. Thus, it could be interesting to minimize the quenching time with the aid of
of controls in certain cases.
It deserves to mention that  the  quenching time $T_q(f,y^0, u)$, with some control $u$,  can really be strictly less than $T_q(f,y^0,0)$. Here is an example. Consider the problem
$(P)^{f^{(2)}}_{y^0}$, where $y_1^0=3/4,\ y_2^0=0$ and $B(t)=
\begin{pmatrix} 1 &
0 \\
0 &
0
\end{pmatrix}
$ for all
$t\in R^+$. It can be directly checked that $T_q(f^{(2)},(3/4,0)^T, (1,0)^T)=1/32<T_q(f^{(2)},(3/4,0)^T, 0)={-1/4}-\ln(3/4)$.

Because of the quenching behavior of solutions to systems (1.1), the
usual methods applied to solve the general time optimal control problems (see, for
instance, \cite{Barbu}, \cite{Fattorini}, \cite{Li-Yong}, \cite{Phung-Wang-Zhang}, \cite{G. Wang}, \cite{J. Yong}
) do not work for
Problem $(P)_{y^0}^{f}$ with $f\in \Lambda$ and $y^0\in S^{f}$.
We approach our main results by the following steps: First, we show some
invariant properties for solutions of systems (\ref{1.1}). Next, we
 prove that given $f\in \Lambda$ and $y^0\in S^{f}$, the corresponding
solution of system (\ref{1.1}), quenches at finite time for each
control $u\in \mathcal{U}_{ad}$. Then we give the quenching rate
estimate for solutions of systems (\ref{1.1}). Furthermore, we
obtain the following property: when a sequence
$\{u_k\}$ of controls tends to a control $u$ in a suitable topology, the
solutions $y(\cdot\ ; f,y^0,u_k)$ with $k$ sufficiently large share a
common interval of non-quenching with the solution $y(\cdot\ ;f, y^0,u)$.
Finally, we use the above-mentioned results to  verify our main theorems.

It deserves to mention the paper \cite{Lin-Wang} where minimal blowup time optimal control problems of
ordinary differential equations were studied and some methods differing from methods to study the usual
time optimal control problems were developed. Though we borrow some ideas from \cite{Lin-Wang}, the current study
 differs from \cite{Lin-Wang} from several points of view as follows:  $(i)$  The study \cite{Lin-Wang} concerns  optimal controls to
minimize the blowup time while our study aims at optimal controls to
minimize the quenching time;  $(ii)$ The methods used to study the problems in \cite{Lin-Wang} do not fit with the problems in the current paper. We develop
new methods in this paper.

The rest of the paper is organized as follows:  Section 2 presents some preliminary lemmas which supply
some properties of solutions to controlled systems (1.1).  Section 3 proves the existence of
optimal controls for Problem $(P)_{y^0}^{f}$ with $f\in \Lambda$ and $y^0\in S^f$. Section 4 verifies Theorem 1.2 and Theorem 1.3.

\section{Preliminary  Lemmas}
\ \ \ \ In this section, we shall introduce some properties of solutions to systems (\ref{1.1}), which will
play important roles to prove our main results.
\subsection{The existence and invariant property }

\ \ \ \ \ Consider the system 
\begin{align}
\label{2.1} \left\{\begin{array}{ll} \displaystyle\frac{d\xi(t)}{dt}
=f(\xi(t))+B(t)u(t),\ \ \ \   t>t_0,\\
\displaystyle \xi(t_0)=y^0.
\end{array}\right.
\end{align}
Since for each $f\in \Lambda$, $f$ is continuously
differential over the domain $S^f$,  it is clear
that given $t_0\geq 0$, $f\in \Lambda$, $y^0\in S^{f}$  and $u\in \mathcal{U}_{ad}$, system (\ref{2.1})
has  a unique solution. We denote this solution  by $\xi(\cdot\ ; t_0,f, y^0,u)
=(\xi_1(\cdot\ ; t_0,f, y^0,u),\ \xi_2(\cdot\ ; t_0,f, y^0,u))^T$, and write $[t_0,T_{max}{(t_0,f,y^0,u)})$ for its maximal interval of
existence.

Let $t_0\geq 0$,  $y^0=(y_1^0,y_2^0)^T\in S^{f^{(1)}}$  and $u\in \mathcal{U}_{ad}$. If $y_1^0<1$, then by the continuity of the solution $\xi(\cdot\ ; t_0,f^{(1)}, y^0,u)$, we can find a
positive number $\alpha$ sufficiently small such that  $\xi_1(t; t_0,f^{(1)}, y^0,u)<1$ for all $t\in [t_0,t_0+\alpha)$.
Similarly,  if $y_1^0>1$, we can find a
positive number $\beta$ sufficiently small such that  $\xi_1(t; t_0,f^{(1)}, y^0,u)>1$ for all $t\in [t_0,t_0+\beta)$.
 Given $t_0\geq 0$, $y^0=(y_1^0,y_2^0)^T\in S^{f^{(1)}}$  and $u\in \mathcal{U}_{ad}$, we shall use the notation $I(t_0,f^{(1)}, y^0,u)$
 to denote such a subinterval of the interval $[t_0,T_{max}{(t_0,f^{(1)},y^0,u)})$: When $y_1^0<1$, $I(t_0,f^{(1)}, y^0,u)$ denotes the
maximal time interval in
which $\xi_1(t; t_0,f^{(1)}, y^0,u)<1$; While  when $y_1^0>1$, $I(t_0,f^{(1)}, y^0,u)$ denotes
the maximal time interval in
which $\xi_1(t; t_0,f^{(1)}, y^0,u)>1$.  By the continuity of the solution $\xi( \cdot\ ; t_0,f^{(1)}, y^0,u)$ again, it is clear from the existence theorem and
the extension theorem of ordinary differential equations that $I(t_0,f^{(1)}, y^0,u)$ is a left closed and right open interval, whose left
end point is $t_0$. Let $1-\frac{1}{2K_0}<K_1<1$,
$1<\widetilde{K}_1<1+\frac{1}{2K_0}$, $K_2>K_0+\frac{1}{K_0}-1$ and
$\widetilde{K}_2>K_0+1$.

We have the the following lemma.\\
{\bf Lemma 2.1.1}
{\it  Given $t_0\geq 0$, $y^0\in S^{f^{(1)}}$ and $u\in \mathcal{U}_{ad}$, it holds that for all $t\in I(t_0,f^{(1)}, y^0,u)$,
\begin{align}
\label{2.2}&\xi_1(t; t_0,f^{(1)}, y^0, u)\geq K_1,\ \xi_2(t; t_0,f^{(1)}, y^0,
u)\geq K_2, \ \mbox{when}\ K_1\leq y^0_1<1\ \mbox{and}\
y^0_2\geq K_2;\\
\label{2.3} &\xi_1(t; t_0, f^{(1)},{y}^0, u)\leq \widetilde{K}_1,\ \xi_2(t;
t_0,f^{(1)},{y}^0, u)\geq \widetilde{K}_2,\ \mbox{when}\ 1<{y}^0_1\leq\widetilde{K}_1\   \mbox{and}\ {y}^0_2\geq \widetilde{K}_2.
\end{align}}

{\bf Proof.}
In order to prove (\ref{2.2}), we claim the following property (A):
{\it Suppose that $t_0\geq 0$,
$y^0\in S^{f^{(1)}}$ with $K_1\leq y^0_1<1$,
$y^0_2\geq K_2$ and $u\in \mathcal{U}_{ad}$. Then, there
is a positive number  $\eta$ with $[t_0, t_0+\eta )\subset I(t_0,f^{(1)}, y^0, u)$ such that the solution $\xi( \cdot\
;t_0,f^{(1)}, y^0, u)$ holds the following inequalities: $\xi_1(t; t_0,f^{(1)}, y^0,
u)\geq K_1$ and $\xi_2(t; t_0,f^{(1)},y^0, u)\geq K_2$ for all $t$ in the
interval $[t_0,t_0+\eta]$.}

The proof of the property is given in what follows. Since $1-\frac{1}{2K_0}<K_1\leq
y^0_1<1$ and $y^0_2\geq K_2>K_0+\frac{1}{K_0}-1$, we can use the continuity of the solution
$\xi(\cdot\ ; t_0, f^{(1)},y^0, u)$ to find a positive constant $\eta$ such
that $[t_0, t_0+\eta )\subset I(t_0,f^{(1)}, y^0, u)$, $1-\frac{1}{2K_0}<\xi_1(t; t_0,f^{(1)},
y^0,u)<1$ and $\xi_2(t; t_0,f^{(1)}, y^0,u)>K_0+\frac{1}{K_0}-1$ for all $t\in [t_0,t_0+\eta]$.
On the other hand, by the definition $K_0$ in (\ref{1.1-1}), it holds that for each $u\in \mathcal{U}_{ad}$,
\begin{align}\label{1.3-1}|b_{1}(t, u(t))|\leq K_0 \mbox{ and}\ |b_{2}(t, u(t))|\leq K_0\ \ \mbox{for a.e.}\ t\in R^+.\end{align}
 Hence, it follows from (\ref{2.1}) with
$f=f^{(1)}$ and the inequality $K_0+\frac{1}{K_0}-1\geq 1$ that for all $t\in [t_0, t_0+\eta]$,
\begin{align}
\xi_1(t)-y^0_1&=\int_{t_0}^t\Big(
\frac{\xi_2(\tau)}{1-\xi_1(\tau)}+b_{1}(\tau, u(\tau))\Big)d\tau\nonumber\\
&\geq\int_{t_0}^t \Big(\frac{K_0+1/K_0-1}{1/(2K_0)}-K_0\Big)d\tau.
\nonumber\\&\geq 0;\nonumber\end{align}
\begin{align}
\xi_2(t)-y^0_2&=\int_{t_0}^t\big(
\xi_1(\tau)+\xi_2(\tau)+b_2(\tau, u(\tau)\big)d\tau\nonumber\\
&\geq\int_{t_0}^t \big(1-1/(2K_0)+K_0+1/K_0-1-K_0\big)d\tau.
\nonumber\\&\geq 0.\nonumber
\end{align}
Here, we simply write $\xi(\cdot)$ for the
solution $\xi(\cdot\ ;t_0,f^{(1)}, y^0, u)$. From the two inequalities
mentioned above and the inequalities $
y^0_1\geq K_1$ and $y^0_2\geq K_2$, we get the property (A).

{\it Now, we come back to prove the property (\ref{2.2}).} 

By seeking a
contradiction, suppose that there existed a time $t_0\geq 0$, an initial data $y^0\in S^{f^{(1)}}$,
a control $u\in \mathcal{U}_{ad}$ and a number pair $(K_1,K_2)$
with $1-{1}/{(2K_0)}<K_1<1$ and $K_2>K_0+1/K_0-1$ such that the
solution $\xi( \cdot\ ;t_0,f^{(1)}, y^0, u)$ with $K_1\leq y^0_1<1$ and
$y^0_2\geq K_2$ does not satisfy (\ref{2.2}). Then we would find a
number $s_0>t_0$ in the interval  $I(t_0,f^{(1)},y^0, u)$ such
that   $\xi_1(s_0;t_0,f^{(1)}, y^0, u)< K_1$ or $\xi_2(s_0;t_0,f^{(1)}, y^0, u)<K_2$. Write $s_1$ for $\inf\{t\in[t_0, s_0 ]; \ \xi_1(t
;t_0,f^{(1)}, y^0, u)< K_1\}$  and write $s_2$ for $\inf\{t\in[t_0, s_0 ]; \
\xi_2(t;t_0,f^{(1)}, y^0, u)< K_2\}.$ We may as well assume that $s_1\leq
s_2$. Since $\xi(\cdot\ ;t_0,f^{(1)}, y^0, u)$ is continuous over $I(t_0,f^{(1)},y^0, u)$, we
can use the definition of $s_1$ and $s_2$, and the inequality
$s_1\leq s_2$ to derive the following properties: (a)
$\xi_1(s_1;t_0,f^{(1)}, y^0, u)\geq K_1$; (b) Corresponding to each $\delta$ with $(s_1,
s_1+\delta )\subset I(t_0,f^{(1)},y^0, u)$, there exists a number $t_\delta$ in $(s_1,
s_1+\delta )$ such that $\xi_1(t_\delta;t_0,f^{(1)}, y^0, u)<K_1$; (c) $\xi_2(t;t_0,f^{(1)}, y^0, u)\geq
K_2$, for each $t\in [t_0,s_1]$. Write $z=\xi(s_1;t_0,f^{(1)}, y^0,u)$.
Then, it follows from the property (a), (c) and the definition of
the interval $I(t_0,f^{(1)},y^0,\\ u)$ that $K_1\leq z_1=\xi_1(s_1;t_0,f^{(1)}, y^0, u)<1$ and $z_2=\xi_2(s_1;t_0,f^{(1)}, y^0, u )\geq K_2$. Consider the
following system,
\begin{align}
\label{2.4-1}   \left\{\begin{array}{ll}
\displaystyle\frac{d\xi(t)}{dt}=f(\xi(t))+B(t)u(t),\ \ \ \   t>s_1,\\
\displaystyle\xi(s_1)=z.\end{array}\right.
\end{align}
Making use of the
property (A), we can choose a positive number $\delta_0$
sufficiently small such that
\begin{align}
\label{2.4-2}\xi_1(t; s_1,f^{(1)}, z, u)\geq K_1,\ \xi_2(t; s_1,f^{(1)}, z, u)\geq
K_2\ \mbox{for each}\ t\in[ s_1, s_1\mbox{$+$}\delta_0 ).
\end{align}
On the other hand,  the function $\xi( \cdot\ ;t_0,f^{(1)}, y^0, u)$ also
satisfies the system (\ref{2.4-1}) when $t\geq s_1$. By the uniqueness of the solution to the system
(\ref{2.4-1}), we necessarily have that $ \xi(t;{t}_0,f^{(1)}, y^0, u)=
\xi(t; s_1,f^{(1)}, z, u )$ for all $t\geq s_1$. This, combined with
the above-mentioned property {(b)}, imply the existence of a number
$t_{\delta_0}$ in $( s_1, s_1+\delta_0)$ such that $
\xi_1(t_{\delta_0}; s_1,f^{(1)}, z,u)<K_1.$ This contradicts to
(\ref{2.4-2}). Therefore, we have accomplished the proof of the
property (\ref{2.2}).

{\it Next, we can use the very similar argument to prove  (\ref{2.3}), we shall
only give a key point of the proof of (\ref{2.3}).}

Write ${\xi}(\cdot)$ for solution $\xi( \cdot\ ; t_0,f^{(1)},{y}^0,u)$ to
the system (\ref{2.1}), where $t_0\geq 0$, $y^0\in S^{f^{(1)}}$ with $1<{y}^0_1\leq\widetilde{K}_1$,
${y}^0_2\geq \widetilde{K}_2$ and $u\in \mathcal{U}_{ad}$. In this case,
we can use the similar argument we used to prove the above-mentioned property
(A) to get the following property: we can find a positive
constant $\widetilde{\eta}$ such that $[t_0, t_0+\widetilde{\eta}
)\subset I(t_0,f^{(1)}, {y}^0, u)$, $1<\xi_1(t)<1+{1}/{(2K_0)}$ and
$\xi_2(t)>K_0+1$ for all $t\in [t_0, t_0+\widetilde{\eta}]$.
Then, from (\ref{2.1}) with $f=f^{(1)}$ and (\ref{1.3-1}), it holds that for all $t\in [t_0, t_0+\widetilde{\eta}]$,
\begin{align}
\xi_1(t)-{y}^0_1&=\int_{t_0}^t\Big(
\frac{{\xi}_2(\tau)}{1-{\xi}_1(\tau)}+b_{1}(\tau, u(\tau)\Big)d\tau
\nonumber\end{align}\begin{align}
&\leq\int_{t_0}^t \Big(-\frac{K_0+1}{1/(2K_0)}+K_0\Big)d\tau. \nonumber\\&\leq
0;\nonumber\end{align}
\begin{align}
\xi_2(t)-y^0_2&=\int_{t_0}^t\big(
{\xi}_1(\tau)+{\xi}_2(\tau)+b_{2}(\tau, u(\tau)\big)d\tau\nonumber\\
&\geq\int_{t_0}^t \big(1+1+K_0-K_0\big)d\tau. \nonumber\\&\geq 0.\nonumber
\end{align} 

Further, we can use the similar argument to prove (\ref{2.2}) to complete the proof of
(\ref{2.3}).
\ \ \ \#\\

Let $t_0\geq 0$,  $y^0\in S^{f^{(2)}}$  and $u\in \mathcal{U}_{ad}$, we shall use the notation $I(t_0,f^{(2)}, y^0,u)$
 to denote such a subinterval of the interval $[t_0,T_{max}{(t_0,f^{(2)},y^0,u)})$: When $\|y^0\|<1$, $I(t_0,f^{(2)}, y^0,u)$ denotes
the maximal time interval in
which $\|\xi(t; t_0,f^{(2)}, y^0,u)\|<1$; While when $\|y^0\|>1$, $I(t_0,f^{(2)}, y^0,u)$ denotes
the maximal time interval in
which $\|\xi(t; t_0,f^{(2)}, y^0,u)\|>1$. By the continuity of the solution $\xi( \cdot\ ; t_0,f^{(2)}, y^0,u)$, it is clear from the existence theorem and
the extension theorem of ordinary differential equations that $I(t_0,f^{(2)}, y^0,u)$ is a left closed and right open interval, whose left
end point is $t_0$. Let $1-\frac{1}{2K_0+1}<K_3<1$ and
$1<\widetilde{K}_3<1+\frac{1}{2K_0}$.

We have the following lemma.\\
{\bf Lemma 2.1.2}\;\;
{\it Given $t_0\geq 0$, $y^0\in S^{f^{(2)}}$ and $u\in \mathcal{U}_{ad}$, it holds that for all $t\in I(t_0,f^{(2)}, y^0,u)$,
\begin{align}
\label{2.6}&\|\xi(t; t_0,f^{(2)}, y^0, u)\|\geq K_3,\ \ \mbox{when}\
K_3\leq \|y^0\|<1;\\
\label{2.7} &\|\xi(t; t_0,f^{(2)},y^0, u)\|\leq \widetilde{K}_3,\ \ \mbox{when}\
1<\|y^0\|\leq \widetilde{K}_3.
\end{align}}

{\bf Sketch proof of Lemma 2.1.2.}

We can use the very similar argument of the proof of Lemma 2.1.1 to prove this lemma. Here, we only give
the key procedures to prove this lemma as follows.

{\it To prove (\ref{2.6}) }

 Write $\xi(\cdot)$ for solution $\xi( \cdot\ ; t_0,f^{(2)},y^0,u)$, where $t_0\geq 0$,
$y^0\in S^{f^{(2)}}$ with $K_3\leq \|y^0\|<1$ and $u\in \mathcal{U}_{ad}$. By system (\ref{2.1}) with $f=f^{(2)}$, it
holds that
\begin{align}
\label{2.8} \left\{\begin{array}{ll} \displaystyle\frac{d\|\xi(t)\|}{dt}
=\displaystyle\frac{\|\xi(t)\|}{1-\|\xi(t)\|}+<\displaystyle\frac{\xi(t)}{\|\xi(t)\|}, B(t)u(t)>,\ \ \ \
t>t_0,\\\displaystyle \|\xi(t_0)\|=\|y^0\|.
\end{array}\right.
\end{align}
The key procedure to prove (\ref{2.6})  is that we can use the similar argument
of the proof of property (A) of Lemma 2.1.1 to get the following property: we can find a positive constant
$\eta$ such that $[t_0, t_0+\eta )\subset I(t_0,f^{(2)}, y^0, u)$ and  $1-\frac{1}{2K_0+1}<\|\xi(t)\|<1$ for all $t\in [t_0,t_0+{\eta}]$.
Then, from (\ref{2.8})
and (\ref{1.3-1}), it holds that for all
$t\in [t_0,t_0+{\eta}]$,
\begin{align}
\|\xi(t)\|-\|y^0\|&=\int_{t_0}^t\Big\{\displaystyle\frac{\|\xi(\tau)\|}{1-\|\xi(\tau)\|}
+<\displaystyle\frac{\xi(\tau)}{\|\xi(\tau)\|}, B(\tau)u(\tau)>\Big\}d\tau
\nonumber\\
&\geq\int_{t_0}^t \Big(\frac{1-\frac{1}{2K_0+1}}{\frac{1}{2K_0+1}}-K_0\Big)d\tau
\nonumber\\&\geq 0.\nonumber\end{align}
Further, we can use the similar argument of the proof of Lemma 2.1.1 to complete the proof of
(\ref{2.6}) of this lemma.

{\it To prove (\ref{2.7}) }

Write ${\xi}(\cdot)$ for solution $\xi(\cdot\ ; t_0,f^{(2)},{y}^0,u)$, where $t_0\geq 0$,
$y^0\in S^{f^{(2)}}$ with $1<\|y^0\|\leq \widetilde{K}_3$ and $u\in \mathcal{U}_{ad}$.
The key procedure to prove (\ref{2.7}) is that we can use the similar argument of the proof of the
property (A) of Lemma 2.1.1 to get the following property: we can find a positive constant $\widetilde{\eta}$
such that $[t_0, t_0+\widetilde{\eta} )\subset I(t_0,f^{(2)}, y^0, u)$,  $1<\|\xi(t)\|<1+\frac{1}{2K_0}$ for all $t\in [t_0,
t_0+\widetilde{\eta}]$.
Then, from (\ref{2.8}) and (\ref{1.3-1}), it holds that for all $t\in [t_0,t_0+\widetilde{\eta}]$,
\begin{align}
\|{\xi}(t)\|-\|{y}^0\|&=\int_{t_0}^t\Big\{\displaystyle\frac{\|{\xi}(\tau)\|}{1-\|{\xi}(\tau)\|}
+<\displaystyle\frac{{\xi}(\tau)}{\|{\xi}(\tau)\|}, B(\tau)u(\tau)>\Big\}d\tau
\nonumber\\&\leq
\int_{t_0}^t (-2K_0+K_0)d\tau
\nonumber\\&\leq 0.\nonumber
\end{align}
Further, we can use the similar argument of the proof of Lemma 2.1.1 to complete the proof of
(\ref{2.7}).\ \ \ \#
\\

Let $t_0\geq 0$, $y^0\in S^{f^{(3)}}$  and $u\in \mathcal{U}_{ad}$, we shall use the notation $I(t_0,f^{(3)}, y^0,u)$
 to denote such a subinterval of the interval $[t_0,T_{max}{(t_0,f^{(3)},y^0,u)})$: When $y_1^0<1$ and $y_2^0<1$, $I(t_0,f^{(3)}, y^0,u)$ denotes
the maximal time interval in
which $\xi_1(t; t_0,f^{(3)}, y^0,u)<1$ and $\xi_2(t; t_0,f^{(3)}, y^0,u)<1$. While when $y_1^0>1$ and $y_2^0>1$, $I(t_0,f^{(3)}, y^0,u)$ denotes
the maximal time interval in
which $\xi_1(t; t_0,f^{(3)}, y^0,u)>1$ and $\xi_2(t; t_0,f^{(3)}, y^0,u)>1$. By the continuity of the solution $\xi( \cdot\ ; t_0,f^{(3)}, y^0,u)$, it is clear from the existence theorem and
the extension theorem of ordinary differential equations that $I(t_0,f^{(3)}, y^0,u)$ is a left closed and right open interval, whose left
end point is $t_0$. Let $1-\frac{e^{-3/2}}{2K_0}<K_4<1$ and $1<\widetilde{K}_4<1+\frac{e^{-3/2}}{2K_0}$.

We have the the following lemma.\\
{\bf Lemma 2.1.3}\ \
{\it  Given $t_0\geq 0$, $y^0\in S^{f^{(3)}}$ and $u\in \mathcal{U}_{ad}$, it holds that for all $t\in I(t_0,f^{(3)}, y^0,u)$,
\begin{align}
\label{2.9}&\xi_1(t; t_0,  f^{(3)},y^0, u)\geq K_4,\ \xi_2(t; t_0, f^{(3)}, y^0, u)\geq K_4,\ \mbox{when}\
K_4\leq y^0_1<1,\ K_4\leq y^0_2<1;
\end{align}
\begin{align}
\label{2.10} &\xi_1(t; t_0, f^{(3)}, {y}^0, u)\leq \widetilde{K}_4,\ \xi_2(t; t_0, f^{(3)}, {y}^0, u)\leq \widetilde{K}_4,
\ \mbox{when}\
1<y^0_1\leq \widetilde{K}_4,\ 1<y^0_2\leq \widetilde{K}_4.\end{align}
}

Indeed, we can use the very similar argument of the proof of Lemma 2.1.1 to prove Lemma 2.1.3.\\

\subsection{Quenching property and estimate of quenching rate}

{\bf Lemma 2.2.1}
{\it For each $y^0\in S^{f^{(1)}}$ and each $u\in \mathcal{U}_{ad}$, there exists a
time $T_q(f^{(1)},y^0,u)$ in the interval $(0, T_{max}{(f^{(1)},y^0,u)}]$ holding the property that \begin{align}\label{2.13-1}
T_q(f^{(1)},&y^0,u)\leq(y_1^0-1)^2<+\infty,\\ \label{2.11}
\ \lim\limits_{t\rightarrow T_q(f^{(1)},y^0,u)}y_1(t;f^{(1)},y^0,u)=1&,\ \lim\limits_{t\rightarrow T_q(f^{(1)},y^0,u)}\|f^{(1)}(y(t;f^{(1)},y^0,u))\|=+\infty,
\\
\label{2.13-1-1}\|f^{(1)}(y(t;f^{(1)},y^0,u))\|&<+\infty\;\;\mbox{as}\;\; t\in [0,T_q(f^{(1)},y^0,u)).\end{align}
  Moreover, there exists a positive constant
$C$, independent of $y^0\in S^{f^{(1)}}$ and $u$,  such that
\begin{align}
\label{2.12}\displaystyle\frac{1}{|1-y_1(t; f^{(1)},y^0, u)|}\leq
C(T_q(f^{(1)},y^0,u)-t)^{-{2}/{3}}\ \ \mbox{for each}\ t\in [0,
T_q(f^{(1)},y^0,u)).\end{align}}

{\bf Proof.}
Suppose that $y^0\in S^{f^{(1)}}$ and $u\in \mathcal{U}_{ad}$. Let $T_q(f^{(1)},y^0,u)$ be the right end point of
the interval $I(0,f^{(1)},y^0,u)$. Thus, $I(0,f^{(1)},y^0,u)=[0,T_q(f^{(1)},y^0,u))$ and $T_q(f^{(1)},y^0,u)\in(0, T_{max}{(f^{(1)},y^0,u)}]$
(see the definition and the property
of $I(0,f^{(1)},y^0,\\u)$ on Page 6).
We shall prove Lemma 2.2.1 in the following two cases.

{\it Case 1. $1-\frac{1}{2K_0}<
y^0_1<1$ and $y^0_2>K_0+\frac{1}{K_0}-1$}

We shall complete the proof of Case 1 by the following
three steps.

{\it Step 1. To prove (\ref{2.13-1}) in Case 1} 

By property (\ref{2.2}) in Lemma 2.1.1 and the definition of $I(0,f^{(1)},y^0,u)$ and $T_q(f^{(1)},\\y^0,u)$,
 the solution $y(\cdot\ ;f^{(1)},y^0,u)$
with $1-\frac{1}{2K_0}<
y^0_1<1$, $y^0_2>K_0+\frac{1}{K_0}-1$ and $u\in \mathcal{U}_{ad}$ holds the property that for each $t\in [0,
T_q(f^{(1)},y^0,u))$,
\begin{align}
\label{2.15-1}
1-\frac{1}{2K_0}<y_1^0\leq y_1(t;f^{(1)},y^0,u)<1\ \mbox{and}\ y_2(t;f^{(1)},y^0,u)\geq y_2^0>K_0+\frac{1}{K_0}-1.
\end{align}
Then, from system (\ref{1.1}) with $f=f^{(1)}$, (\ref{1.3-1}) and the inequality $K_0+1/K_0-1\geq 1$, it holds that
\begin{align}
\label{2.15}\displaystyle\frac{dy_1(t)}{dt}&=\displaystyle\frac{y_2(t)}{2(1-y_1(t))}+ \displaystyle\frac{y_2(t)}{2(1-y_1(t))}+b_1(t,u(t))\nonumber\\&\geq\displaystyle\frac{1}
{2(1-y_1(t))}+K_0-K_0\nonumber\\&= \displaystyle\frac{1}
{2(1-y_1(t))}\ \ \ \ \ \mbox{for a.e.}\ t\in [0, T_q(f^{(1)},{y}^0,u)).
\end{align}
Here and throughout the proof, we simply write $y(\cdot)$ for $y(\cdot\ ;f^{(1)},y^0,u)$.

Let $\chi(\cdot)$ be the solution to the following equation:
\begin{align}
\label{2.17}
\left\{\begin{array}{ll}\displaystyle\frac{d\chi(t)}{dt}
=\frac{1}{2(1-\chi(t))},\ \ \ \   t>0,\\
\chi(0)={y}^0_1.\end{array}\right.
\end{align}
Then, we can directly solve the equation (\ref{2.17}) to get that
$$\chi(t)=1-\sqrt{(y_1^0-1)^2-t},\ t\in [0, (y^0_1-1)^2),$$
from which it is easy to check that
\begin{align}
\label{2.17-1}\chi(t)<1,\  t\in [0, (y^0_1-1)^2)\ \ \mbox{and}\ \ \chi(\cdot)\ \mbox{quenches}\ \mbox{at the time}\ (y^0_1-1)^2.\end{align}
Furthermore, making use of (\ref{2.15}) and (\ref{2.17}), we can derive, from
the comparison theorem of ordinary differential equations,
the following inequality:
\begin{equation}y_1(t)\geq \chi(t),\;\; t\in [0,  \mbox{min}\{T_q(f^{(1)},{y}^0, u),
(y^0_1-1)^2\ \} ).\nonumber\end{equation} This, combined with (\ref{2.15-1}) and (\ref{2.17-1}),
implies that (\ref{2.13-1}) in Case 1.

{\it Step 2. To prove (\ref{2.11}) in Case 1}

Indeed, by (\ref{2.15-1}) and (\ref{2.15}), it is clear that $y_1(\cdot)$
is monotonously increasing over the interval $[0, T_q(f^{(1)},y^0,
u))$. This, combined with (\ref{2.15-1}) implies that
$\lim\limits_{t\rightarrow T_q(f^{(1)},y^0, u)} y_1(t)$
exists and $\lim\limits_{t\rightarrow T_q(f^{(1)},y^0, u)}
y_1(t)\leq 1$.

Now, we claim that
\begin{align}\label{2.19}\lim\limits_{t\rightarrow T_q(f^{(1)},y^0, u)}
y_1(t)=1.
\end{align}

By contradiction, if $\lim\limits_{t\rightarrow T_q(f^{(1)},y^0, u)} y_1(t)=\beta<1$. Then, by the
continuity and the monotonicity of the solution  $y_1(\cdot)$, it holds that $y_1(t)\leq\beta$
for $t\in [0,T_q(f^{(1)},y^0, u))$ and $y_1(T_q(f^{(1)},y^0, u))=\beta<1$. Then, we can extend the solution
$y(\cdot)$ and can find an interval $[T_q(f^{(1)},y^0, u), T_q(f^{(1)},y^0, u) +\delta_1)$ with
$\delta_1$ sufficiently small such that $y_1(t)<1$ on the interval $[T_q(f^{(1)},y^0, u), T_q(f^{(1)},y^0, u)
+\delta_1)$. However, $I(0,f^{(1)},y^0,u)=[0,T_q(f^{(1)},y^0,u))$ and $I(0,f^{(1)},y^0,u)$
is the maximal interval in which $y_1(t)<1$. This is a contradiction. Thus, (\ref{2.19}) holds.

On the other hand, by (\ref{2.15-1}), it is clear that
$$\|f^{(1)}(y(t))\|\geq |\frac{y_2(t)}{1-y_1(t)}|> \frac{K_0+\frac{1}{K_0}-1}{1-y_1(t)}\geq \frac{1}{1-y_1(t)},\ \ t\in [0,T_q(f^{(1)},y^0,u)),$$
from which and (\ref{2.19}), it holds that $\lim\limits_{t\rightarrow  T_q(f^{(1)},y^0,u)}
\|f^{(1)}(y(t))\|=+\infty$. This completes the proof of (\ref{2.11}) in Case 1.

{\it Step 3. To prove (\ref{2.13-1-1}) in Case 1}

Since $y_1(t)<1$ for each $t\in [0,T_q(f^{(1)},y^0,u))$ and $y(\cdot)$ is continuous over
the interval $[0,T_q(f^{(1)},y^0,u))$, we can get $(\ref{2.13-1-1})$ in Case 1.

{\it Step 4. To prove (\ref{2.12}) in Case 1}

By system (\ref{1.1}) with
$f=f^{(1)}$, we conclude from (\ref{2.11}), (\ref{2.15-1})
 and  (\ref{1.3-1}) that for each
$t\in [0, T_q(f^{(1)},{y}^0,u))$,
\begin{align}
&-\frac{2}{3}\int_t^{T_q(f^{(1)},y^0,u)}d(
1-y_1(\tau))^{{3}/{2}}\nonumber\\=&\frac{2}{3}(1-y_1(t))^{3/2}\nonumber\\=&\int_t^{T_q(f^{(1)},y^0,u)} (1-y_1(\tau))^{1/2}\Big(\displaystyle\frac{y_2(\tau)}{1-y_1(\tau)}+b_{1}(\tau,u(\tau))\Big)d\tau
\nonumber
\\= & \int_t^{T_q(f^{(1)},y^0,u)}\Big(\displaystyle\frac{y_2(\tau)}{(1-y_1(\tau))^{1/2}}+(1-y_1
(\tau))^{1/2}b_{1}(\tau,u(\tau))\Big)d\tau\nonumber
\\ \geq & \int_t^{T_q(f^{(1)},y^0,u)}\Big(\displaystyle\frac{1}{(1-1+\frac{1}{2K_0})^{{1/2}}}-\frac{K_0}{\sqrt{2K_0}}\Big)d\tau
\nonumber\\ \geq & (\sqrt{2}-\frac{1}{\sqrt{2}})\sqrt{K_0}(T_q(f^{(1)},y^0,u)-t)\nonumber
,\end{align}
which implies (\ref{2.12}) in Case 1.

{\it Case 2. $y^0\in S^{f^{(1)}}$ with $1<
y^0_1<1+\frac{1}{2K_0}\ \mbox{and}\ y^0_2>K_0+1$}

Now we shall give a sketch proof of Case 2, which is similar to the
proof of Case 1 of this lemma.

{\it Step 1. To prove (\ref{2.13-1}) in Case 2}

By property (\ref{2.3}) in Lemma 2.1.1 and the definition of $I(0,f^{(1)},y^0,u)$ and $T_q(f^{(1)},\\y^0,u)$, the solution $y(\cdot\ ;f^{(1)},y^0,u)$
with $1<
y^0_1<1+\frac{1}{2K_0},\ y^0_2>K_0+1$ and $u\in \mathcal{U}_{ad}$ holds the property that for each $t\in [0,
T_q(f^{(1)},y^0,u))$,
\begin{align}
\label{2.20}
1<y_1(t;f^{(1)},y^0,u)<1+\frac{1}{2K_0}\ \mbox{and}\ y_2(t;f^{(1)},y^0,u)>K_0+1.
\end{align}
Write $y(\cdot)$ simply for  $y(\cdot\ ;f^{(1)},y^0,u)$.
Then, from system (\ref{1.1}) with $f=f^{(1)}$ and (\ref{1.3-1}), it holds that
\begin{align}
\label{2.21}\displaystyle\frac{d{y}_1(t)}{dt}&=\displaystyle\frac{{y}_2(t)}{2(1-{y}_1(t))}
+ \displaystyle\frac{{y}_2(t)}{2(1-{y}_1(t))}+b_1(t,u(t))
\nonumber\\&\leq \displaystyle\frac{1}{2(1-{y}_1(t))}-K_0+K_0
\nonumber\\&=\displaystyle\frac{1}{2(1-{y}_1(t))}\ \ \ \ \ \mbox{for a.e. }
t\in [0, T_q(f^{(1)},y^0,u)).
\end{align}

Let ${\widetilde{\chi}}(\cdot)$ (with $1<y^0_1<1+\frac{1}{2K_0}$ in this case) be the solution to the equation (\ref{2.17}).
Then, we can directly solve the equation (\ref{2.17}) to get
$$\widetilde{\chi}(t)=1+\sqrt{(y_1^0-1)^2-t},\ t\in [0, (y^0_1-1)^2),$$
from which it is easy to check that
\begin{align}
\label{2.21-1}\widetilde{\chi}(t)>1,\  t\in [0, (y^0_1-1)^2)\ \ \mbox{and}\ \ \widetilde{\chi}(\cdot)\ \mbox{quenches}\ \mbox{at the time}\ (y^0_1-1)^2.\end{align}
Furthermore, making use of (\ref{2.21}) and (\ref{2.17}), we can derive, from
the comparison theorem of ordinary differential equations,
the following inequality:
\begin{equation}y_1(t)\leq\widetilde{\chi}(t),\;\; t\in [0,  \mbox{min}\{T_q(f^{(1)},{y}^0, u),
(y^0_1-1)^2\ \} ).\nonumber\end{equation} This, together with (\ref{2.20}) and
(\ref{2.21-1}), implies (\ref{2.13-1}) in Case 2.

{\it Step 2 We can use the similar argument we used to prove Step 2 and Step 3 in Case 1 of this lemma to prove
 (\ref{2.11}) and  (\ref{2.13-1-1}) hold in Case 2.}

{\it Step 3 To prove (\ref{2.12}) in Case 2}

By system (\ref{1.1}) with $f=f^{(1)}$, we conclude from (\ref{2.11}), (\ref{2.20})
 and (\ref{1.3-1}) that for each
$t\in [0, T_q(f^{(1)},{{y}}^0,u))$,
\begin{align}
&-\frac{2}{3}\int_t^{T_q(f^{(1)},y^0,u)}d(
y_1(\tau)-1)^{{3}/{2}}\nonumber\\=&\frac{2}{3}({y}_1(t)-1)^{3/2}\nonumber\\=&\int_t^{T_q(f^{(1)},{y}^0,u)} ({y}_1(\tau)-1)^{1/2}
\Big(\displaystyle\frac{{y}_2(\tau)}{{y}_1(\tau)-1}-b_{1}(\tau,u(\tau)\Big)d\tau
\nonumber\\=&\int_t^{T_q(f^{(1)},y^0,u)}\Big(\displaystyle\frac{y_2(\tau)}{(y_1(\tau)-1)^{1/2}}
-(y_1(\tau)-1)^{1/2}b_{1}(\tau,u(\tau))\Big)d\tau\nonumber\\\geq & \int_t^{T_q(f^{(1)}y^0,u)}\Big(\displaystyle\frac{1}{(1+\frac{1}{2K_0}-1)^{1/2}}-\frac{K_0}{\sqrt{2K_0}}\Big)d\tau
\nonumber\\\geq  &(\sqrt{2}-\frac{1}{\sqrt{2}})\sqrt{K_0}(T_q(f^{(1)},{y}^0,u)-t)
,\nonumber
\end{align}
which implies (\ref{2.12}) in Case 2.

This completes the proof of Lemma 2.2.1.\ \ \ \#\\
{\it {\bf Remark 2.2.1} Let $y^0\in S^{f^{(1)}}$ and $u\in \mathcal{U}_{ad}$. Since $I(0,f^{(1)},y^0,u)=[0,T_q(f^{(1)},y^0,u))$, we can conclude from the definition of $I(0,f^{(1)},y^0,u)$ that for each $t\in [0,T_q(f^{(1)},y^0,u))$,
\begin{align}
y_1(t;f^{(1)},y^0,u)<1,\ \mbox{when}\ y_1^0<1;\nonumber\\
y_1(t;f^{(1)},y^0,u)>1,\ \mbox{when}\ y_1^0>1.\nonumber
\end{align}}\\
{\it {\bf Lemma 2.2.2}
For each $y^0\in S^{f^{(2)}}$ and each $u\in \mathcal{U}_{ad}$, there exists a
time $T_q(f^{(2)},y^0,u)$ in the interval $(0, T_{max}{(f^{(2)},y^0,u)}]$ holding the property that
 \begin{align}\label{2.22-1}
T_q(f^{(2)},&y^0,u)\leq\frac{(2K_0+1)(\|y^0\|-1)^2}{2K_0}<+\infty,\\ \label{2.22}
\ \lim\limits_{t\rightarrow  T_q(f^{(2)},y^0,u)}\|y(t;f^{(2)}, &y^0, u)\|=1,\ \lim\limits_{t\rightarrow T_q(f^{(2)},y^0,u)}\|f^{(2)}(y(t;f^{(2)},y^0,u))\|=+\infty,
\\ \label{2.22-1-1}
\|f^{(2)}(y(t;&f^{(2)},y^0,u))\|<+\infty\;\;\mbox{as}\;\; t\in [0,T_q(f^{(2)},y^0,u)).\end{align}
Moreover, there exists a positive constant
$C$, independent of $y^0\in S^{f^{(2)}}$ and $u$,  such that
\begin{align}
\label{2.23}\displaystyle\frac{1}{\big|1-\|y(t;f^{(2)}, y^0, u)\|\big|}\leq C(T_q(f^{(2)},y^0,u)-t)^{-{2}/{3}}\ \ \mbox{for each}\
t\in [0, T_q(f^{(2)},y^0,u)).
\end{align}
}

{\bf Proof.} Suppose that $y^0\in S^{f^{(2)}}$ and $u\in \mathcal{U}_{ad}$. Let $T_q(f^{(2)},y^0,u)$ be the right end
point of the interval
$I(0,f^{(2)},y^0,u)$. Thus, $I(0,f^{(2)},y^0,u)=[0,T_q(f^{(2)},y^0,u))$
and $T_q(f^{(2)},y^0,u)\in(0, T_{max}{(f^{(2)},y^0,u)}]$ (see the definition and the property
of $I(0,f^{(2)},y^0,\\u)$ on Page 8).
We shall prove Lemma 2.2.2 in the following two cases.

{\it Case 1. $y^0\in S^{f^{(2)}}$ with $1-\frac{1}{2K_0+1}<\|y^0\|<1$}

We shall prove Case 1 in the following three steps.

{\it Step 1. To prove (\ref{2.22-1}) in Case 1}

Let $y^0\in S^{f^{(2)}}$ with $1-\frac{1}{2K_0+1}<\|y^0\|<1$ and $u\in \mathcal{U}_{ad}$. Write $y(\cdot)$ simply for
$y(\cdot\ ;f^{(2)},y^0,u)$. By (\ref{2.6}) in Lemma 2.1.2 and the definition of $I(0,f^{(2)},y^0,u)$ and $T_q(f^{(2)},y^0,u)$, the solution $y(\cdot)$
satisfies
\begin{align}
\label{2.26}1-\frac{1}{2K_0+1}<\|y(t)\|<1\ \mbox{for each}\ t\in [0, T_q(f^{(2)},y^0,u)).
\end{align}
By system (\ref{1.1}) with $f=f^{(2)}$, it holds that
\begin{align}
\label{2.27} \left\{\begin{array}{ll} \displaystyle\frac{d\|y(t)\|}{dt}
=\displaystyle\frac{\|y(t)\|}{1-\|y(t)\|}+<\displaystyle\frac{y(t)}{
\|y(t)\|}, B(t)u(t)>,\ \ \ \   t>0,\\
\displaystyle \|y(0)\|=\|{y}^0\|,
\end{array}\right.
\end{align}
from which, (\ref{2.26}) and (\ref{1.3-1}), we derive that
\begin{align}\label{2.28}
\displaystyle\frac{d\|y(t)\|}{dt}&\geq \displaystyle\frac{\|y(t)\|}{2(1-\|y(t)\|)}+\displaystyle\frac{\|y(t)\|}
{2(1-\|y(t)\|)}-K_0\nonumber\\&\geq \displaystyle\frac{\|y(t)\|}{2(1-\|y(t)\|)}+ \displaystyle\frac{1-\frac{1}
{2K_0+1}}{2(1-1+\frac{1}{2K_0+1})}-K_0\nonumber\\&\geq \displaystyle\frac{1-\frac{1}{2K_0+1}}{2(1-\|y(t)\|)}\ \
\ \ \ \mbox{for a.e.}\ t\in [0, T_q(f^{(2)},{y}^0,u)).\end{align}

Let $\chi(\cdot)$ be the solution to the following equation:
\begin{align}
\label{2.29}
\left\{\begin{array}{ll}\displaystyle\frac{d\chi(t)}{dt}
=\frac{1-\frac{1}{2K_0+1}}{2(1-\chi(t))},\ \ \ \   t>0,\\
\chi(0)=\|{y}^0\|.\end{array}\right.
\end{align}
Then, we can directly solve the equation (\ref{2.29}) to get  that
$$\chi(t)=1-\sqrt{(\|y^0\|-1)^2-(1-\frac{1}{2K_0+1})t},\ t\in \big[0,\frac{(2K_0+1)
(\|y^0\|-1)^2}{2K_0}\big),$$ from which it is easy to
check that
\begin{align}
\label{2.29-1}\left\{\begin{array}{ll}\chi(t)<1, \ t\in \big[0, \displaystyle{\frac{(2K_0+1)(\|y^0\|-1)^2}{2K_0}}\big)\ \\ \mbox
{and} \ \chi(\cdot)\ \mbox{quenches at the time}\ \displaystyle{\frac{(2K_0+1)(\|y^0\|-1)^2}{2K_0}}.\end{array}\right.\end{align}
Furthermore, making use of (\ref{2.28}) and (\ref{2.29}), we can derive,
from the comparison theorem of ordinary differential equations,
the following inequality:
\begin{equation}\|y(t)\|\geq \chi(t),\;\; t\in \big[0,  \mbox{min}\big\{T_q(f^{(2)},{y}^0, u),
\frac{(2K_0+1)(\|y^0\|-1)^2}{2K_0}\ \big\} \big).\nonumber
\end{equation} This, combined with (\ref{2.26}) and (\ref{2.29-1}), implies (\ref{2.22-1}) in Case 1.

{\it Step 2  We can use the similar argument we used to prove Step 2 and Step 3 in Case 1 of Lemma 2.2.1 to prove
 (\ref{2.22}) and  (\ref{2.22-1-1}) hold in Case 1.}

{\it Step 3. To prove (\ref{2.23}) in Case 1}

By system (\ref{1.1}) with $f=f^{(2)}$, we conclude from (\ref{2.22}), (\ref{2.26}) and (\ref{1.3-1}) that for each
$t\in [0, T_q(f^{(2)},{y}^0,u))$,
\begin{align}
&-\frac{2}{3}\int_t^{T_q(f^{(2)},y^0,u)}d(
1-\|y(\tau)\|)^{{3}/{2}}\nonumber\\=&\frac{2}{3}(1-\|y(t)\|)^{3/2}\ \ \ \ \ \ \ \ \ \ \ \ \nonumber\\=&\int_t^{T_q(f^{(2)},y^0,u)} (1-\|y(\tau)\|)^{1/2}\Big(\frac{\|y(\tau)\|}{1-\|y(\tau)\|}+<\frac{y(\tau)}{\|y(\tau)\|},B(\tau)u(\tau)>\Big)d\tau
\nonumber\\=& \int_t^{T_q(f^{(2)},y^0,u)}\Big(\displaystyle\frac{\|y(\tau)\|}{(1-\|y(\tau)\|)^{{1/2}}}
+(1-\|y(\tau)\|)^{1/2}<\frac{y(\tau)}{\|y(\tau)\|},B(\tau)u(\tau)>\Big)d\tau
\nonumber\\\geq & \int_t^{T_q(f^{(2)},y^0,u)}\Big(\displaystyle\frac{1-\frac{1}{2K_0+1}}{(1-1+\frac{1}{2K_0+1})^{{1/2}}}
-\frac{K_0}{\sqrt{2K_0+1}}\Big)d\tau\nonumber\\
=&\frac{K_0}{\sqrt{2K_0+1}}(T_q(f^{(2)},y^0,u)-t),\nonumber
\end{align}
which implies (\ref{2.23}) in Case 1.

{\it Case 2. $y^0\in S^{f^{(2)}}$ with $1<\|y^0\|<1+\frac{1}{2K_0}$}

The proof of Case 2 is similar to Case 1 of this lemma. Here, we shall give a sketch proof of Case 2.

Let ${y}^0\in S^{f^{(2)}}$ with $1<\|{y}^0\|<1+\frac{1}{2K_0}$
and $u\in \mathcal{U}_{ad}$. Write ${y}(\cdot)$ simply for $y(\cdot\ ;f^{(2)},{y}^0,u)$.
By (\ref{2.7}) in Lemma 2.1.2 and the definition of $I(f^{(2)},{y}^0,u)$ and $T_q(f^{(2)},{y}^0,u)$, the solution ${y}(\cdot)$
satisfies
\begin{align}
\label{2.30}1<\|{y}(t)\|<1+\frac{1}{2K_0}\ \mbox{for each}\ t\in [0, T_q(f^{(2)},{y}^0,u)),
\end{align}
from which, (\ref{2.27}) and (\ref{1.3-1}), we derive that\newpage
\begin{align}
\displaystyle\frac{d\|{y}(t)\|}{dt}&\leq
\displaystyle\frac{\|{y}(t)\|}{2(1- \|{y}(t)\|)}+
\displaystyle\frac{\|{y}(t)\|}{2(1-\|{y}(t)\|)}+K_0 \nonumber\\&\leq
\displaystyle\frac{1}{2(1-\|{y}(t)\|)}\ \ \ \ \ \mbox{for a.e.}\
t\in [0, T_q(f^{(2)},{{y}}^0,u)).\nonumber\end{align}
Then, we can use the similar argument we used to prove Step 1, Step 2 and Step 3 in Case 1 of Lemma 2.2.1 to prove
 (\ref{2.22-1}), (\ref{2.22}) and  (\ref{2.22-1-1}) hold in Case 2 of this lemma.

Now we shall prove (\ref{2.23}) in Case 2.

By system (\ref{1.1}) with $f=f^{(2)}$, we conclude from (\ref{2.22}), (\ref{2.30}) and (\ref{1.3-1}) that
\begin{align}
&-\frac{2}{3}\int_t^{T_q(f^{(2)},{y}^0,u)}d(
\|{y}(\tau)\|-1)^{{3}/{2}}\nonumber\\=&\frac{2}{3}(\|{y}(t)\|-1)^{3/2}
\nonumber\\=&\int_t^{T_q(f^{(2)},{y}^0,u)} (\|{y}(\tau)\|-1)^{1/2}\
\Big(\frac{\|{y}(\tau)\|}{\|{y}(\tau)\|-1}-<\frac{{y}(\tau)}{\|{y}(\tau)\|},B(\tau)u(\tau)>\Big)d\tau
\nonumber\\\geq & \int_t^{T_q(f^{(2)},{y}^0,u)}\Big(\displaystyle\frac{\|{y}(\tau)\|}{(\|{y}(\tau)\|-1)^{{1/2}}}
-\frac{K_0}{\sqrt{2K_0}}\Big)d\tau
\nonumber\\
\geq &(\sqrt{2}-\frac{1}{\sqrt{2}})\sqrt{K_0}(T_q(f^{(2)},{y}^0,u)-t)\ \ \ \ \ \mbox{for each}\
t\in [0, T_q(f^{(2)},{{y}}^0,u)),\nonumber
\end{align}
which implies (\ref{2.23}) in Case 2.\ \ \ \#\\
{\it {\bf Remark 2.2.2} Let $y^0\in S^{f^{(2)}}$ and $u\in \mathcal{U}_{ad}$. Since $I(0,f^{(2)},y^0,u)=[0,T_q(f^{(2)},y^0,u))$, we can conclude from the definition of $I(0,f^{(2)},y^0,u)$ that for each $t\in [0,T_q(f^{(2)},y^0,u))$,
\begin{align}
\|y(t;f^{(2)},y^0,u)\|<1,\ \mbox{when}\ \|y^0\|<1;\nonumber\\
\|y(t;f^{(2)},y^0,u)\|>1,\ \mbox{when}\ \|y^0\|>1.\nonumber
\end{align}}\\
{\it{\bf Lemma 2.2.3}
For each $y^0\in S^{f^{(3)}}$ and each $u\in \mathcal{U}_{ad}$, there exists a
time $T_q(f^{(3)},y^0,u)$ in the interval $(0, T_{max}{(f^{(3)},y^0,u)}]$ holding the property that
\begin{align}\label{2.41}T_q(f^{(3)},y^0,u)\leq\frac{1}{4K_0^2}<+\infty,\ \ \ \ \ \ \ \ \ \ \ \ \ \ \ \ \ \ \ \ \\
\label{2.33}\left\{\begin{array}{ll}\lim\limits_{t\rightarrow  T_q(f^{(3)},y^0,u)}y_1(t;f^{(3)}, y^0, u)=
\lim\limits_{t\rightarrow  T_q(f^{(3)},y^0,u)}y_2(t; f^{(3)},y^0, u)=1,\\
 \lim\limits_{t\rightarrow  T_q(f^{(3)},y^0,u)}\|f^{(3)}(y(t;f^{(3)}, y^0, u))\|=+\infty,\end{array}\right.\\
\|f^{(3)}(y(t;f^{(3)},y^0,u))\|<+\infty\;\;\mbox{as}\;\; t\in [0,T_q(f^{(3)},y^0,u)).\label{2.41-1-1}\end{align}
Moreover, there exists a positive constant
$C$ depending on $y^0\in S^{f^{(3)}}$, but independent of $u$,  such that for each
$t\in [0, T_q(f^{(3)},y^0,u))$,
\begin{align}
\label{2.34}\max\Big\{\displaystyle\frac{1}{|1-y_1(t; f^{(3)},y^0, u)|}, \displaystyle\frac{1}{|1-y_2(t; f^{(3)},y^0, u)|}\Big\}
\leq C(T_q(f^{(3)},y^0,u)-t)^{-{2}/{3}}.
\end{align}
}

{\bf Proof.}  Suppose that $y^0\in S^{f^{(3)}}$ and $u\in \mathcal{U}_{ad}$. Let $T_q(f^{(3)},y^0,u)$ be the right end point of
$I(0,f^{(3)},y^0,u)$. Thus, $I(0,f^{(3)},y^0,u)=[0,T_q(f^{(3)},y^0,u))$ and $T_q(f^{(3)},y^0,u)\in(0, T_{max}{(f^{(3)},y^0,u)}]$ (see the definition and the property
of $I(0,f^{(3)},y^0,u)$ on Page 9). We shall prove Lemma 2.2.3 in the following two cases.

{\it Case 1. $y^0\in S^{f^{(3)}}$ with $1-\frac{e^{-3/2}}{2K_0}<y^0_1<1$ and $1-\frac{e^{-3/2}}{2K_0}<y^0_2<1$}

We shall complete the proof of Case 1 by the following
five steps.

{\it Step 1. To prove (\ref{2.41}) in Case 1}

Let $y^0\in S^{f^{(3)}}$ with $1-\frac{e^{-3/2}}{2K_0}<y^0_1<1$, $1-\frac{e^{-3/2}}{2K_0}<y^0_2<1$ and $u\in \mathcal{U}_{ad}$. Write $y(\cdot)$ simply for $y(\cdot\ ;f^{(3)},y^0,u)$.

By (\ref{2.9}) in Lemma 2.1.3 and the definition of $I(0,f^{(3)},y^0,u)$ and $T_q(f^{(3)},y^0,u)$, the solution $y(\cdot)$
 satisfies
\begin{align}
\label{2.37}
1-\frac{e^{-3/2}}{2K_0}<
y_1(t)<1\ \mbox{and}\ 1-\frac{e^{-3/2}}{2K_0}<
y_2(t)<1\ \mbox{for each}\ t\in [0,
T_q(f^{(3)},y^0,u)).
\end{align}
Then, from system (\ref{1.1}) with $f=f^{(3)}$, the
inequalities $1-\frac{1}{2K_0}<1-\frac{e^{-3/2}}{2K_0}<1$ and (\ref{1.3-1}), it holds that
\begin{align}
\label{2.38}\displaystyle\frac{dy_1(t)}{dt}&= \displaystyle\frac{1}{2(1-y_2(t))}+\displaystyle\frac{1}{2(1-y_2(t))}+b_1(t,u(t))
\nonumber\\&\geq \displaystyle\frac{1}{2(1-y_2(t))}\ \ \ \ \ \mbox{for a.e.}\
t\in [0, T_q(f^{(3)},{y}^0,u)).
\end{align}
Similarly, we have
\begin{align}
\label{2.39}\displaystyle\frac{dy_2(t)}{dt}\geq \displaystyle\frac{1}{2(1-y_1(t))}\ \ \ \ \ \mbox{for a.e.}\
t\in [0, T_q(f^{(3)},{y}^0,u)).
\end{align}
Let $\chi(\cdot)=(\chi_1(\cdot),\chi_2(\cdot))$ be the solution to the following equations:
\begin{align}
\label{2.40}
\left\{\begin{array}{ll}\displaystyle\frac{d\chi_1(t)}{dt}
=\frac{1}{2(1-\chi_2(t))},\ \ \ \   t>0,\\
\displaystyle\frac{d\chi_2(t)}{dt}
=\frac{1}{2(1-\chi_1(t))},\ \ \ \   t>0,\\
\chi_1(0)=\chi_2(0)=1-\frac{1}{2K_0}.\end{array}\right.
\end{align}
Then, it is obvious that
$$\chi_1(t)=\chi_2(t)=1-\sqrt{\frac{1}{4K_0^2}-t},\ t\in \big[0, \frac{1}{4K_0^2}\big),$$
from which it is clear that
\begin{align}
\label{2.40-1}\chi_1(t)<1,\ \chi_2(t)<1, \  t\in \big[0,  \frac{1}{4K_0^2}\big),\ \mbox{both}\ \chi_1(\cdot)
\ \mbox{and}\ \chi_2(\cdot)\ \mbox{quench at the time}\ \frac{1}{4K_0^2}.\end{align}
Furthermore, making use of (\ref{2.38}), (\ref{2.39}) and (\ref{2.40}), we can derive, from
the comparison theorem of ordinary differential equations,
the following inequality:
\begin{equation}y_1(t)\geq \chi_1(t),\;\;y_2(t)\geq \chi_2(t),\ \  t\in \big[0,  \mbox{min}\big\{T_q(f^{(3)},{y}^0, u),
\frac{1}{4K_0^2}\big\} \big).\nonumber\end{equation}
This, together with (\ref{2.37}) and (\ref{2.40-1}), implies that  (\ref{2.41}) in Case 1.

{\it Step 2\ \ We shall prove that $\lim\limits_{t\rightarrow  T_q(f^{(3)},{y}^0,u)}y_1(t)=1$
if and only if $\lim\limits_{t\rightarrow  T_q(f^{(3)},{y}^0,u)}y_2(t)\\=1$.}

First, we shall prove the following
two inequalities.
\begin{align}\label{2.34-1}
\frac{1-y_2(t)}{1-y_1(t)}\leq\frac{1-y^0_2}{1-y^0_1}e^{3/2}\ \ \mbox{for each}\ t\in [0, T_q(f^{(3)},{y}^0, u))
\end{align}
and
\begin{align}\label{2.34-2}
\frac{1-y_1(t)}{1-y_2(t)}\leq\frac{1-y^0_1}{1-y^0_2}e^{3/2}\ \ \mbox{for each} \ t\in [0, T_q(f^{(3)},{y}^0, u)).
\end{align}
Indeed, by system (\ref{1.1}) with $f=f^{(3)}$, it is obvious that
\begin{align}\label{2.42}\left\{\begin{array}{ll} \displaystyle\frac{dy_1(t)}{dt}=
\frac{1}{1-y_2(t)}+b_{1}(t,u(t)),\ \ \ \   t>0,\\
\displaystyle\frac{dy_2(t)}{dt}= \frac{1}{1-y_1(t)}+b_{2}(t,u(t)),\ \ \ \   t>0,\\
 y_1(0)=y_1^0,\ y_2(t)=y_2^0.\end{array}\right.
\end{align}
Thus, we have
\begin{align}
&y_1(t)-y_1^0=\int_0^t \frac{1}{1-y_2(\tau)}d\tau+\int_0^t b_{1}(\tau,u(\tau))d\tau,\ \ t\in
[0, T_q(f^{(3)},{y}^0, u)),\nonumber\\&
y_2(t)-y_2^0=\int_0^t \frac{1}{1-y_1(\tau)}d\tau+\int_0^t b_{2}(\tau,u(\tau))d\tau,\ \ t\in
[0, T_q(f^{(3)},{y}^0, u)).\nonumber
\end{align}
Since $1-\frac{1}{2K_0}<1-\frac{e^{-3/2}}{2K_0}<y_1^0<1$ and $1-\frac{1}{2K_0}<1-\frac{e^{-3/2}}{2K_0}<y_2^0<1$, it holds from the
above two equations, (\ref{2.37}), (\ref{2.41}) and (\ref{1.3-1}) that for each $t\in [0, T_q(f^{(3)},{y}^0, u))$,
\begin{align}\label{2.43}
\int_0^t \frac{1}{1-y_1(\tau)}d\tau=&y_2(t)-y_2^0-\int_0^t b_{2}(\tau,u(\tau))d\tau\nonumber\\
\leq &1-(1-\frac{1}{2K_0})+K_0\cdot \frac{1}{4K_0^2}\nonumber\\\leq &\frac{3}{4K_0}.\end{align}
Similarly, we have
\begin{align}\label{2.431}
\int_0^t \frac{1}{1-y_2(\tau)}d\tau\leq\frac{3}{4K_0}.
\end{align}
On the other hand, for each $t\in [0, T_q(f^{(3)},{y}^0, u))$, multiplying
the first equation of (\ref{2.42}) by $1/(1-y_1(t))$ and integrating
over $[0,t)$, we obtain that
\begin{align}
&-\ln (1-y_1(t))+\ln (1-y^0_1)\nonumber\\=&\int_0^t
\frac{1}{(1-y_1(\tau))(1-y_2(\tau))}d\tau+\int_0^t
\frac{b_{1}(\tau,u(\tau))}{1-y_1(\tau)}d\tau. \nonumber
\end{align}
Similarly, from the second equation of (\ref{2.42}), we get
\begin{align}
&-\ln (1-y_2(t))+\ln (1-y^0_2)\nonumber\\=&\int_0^t
\frac{1}{(1-y_1(\tau))(1-y_2(\tau))}d\tau+\int_0^t
\frac{b_{2}(\tau,u(\tau))}{1-y_2(\tau)}d\tau\ \ \ \mbox{for each}\
t\in [0, T_q(f^{(3)},y^0, u)).\nonumber
\end{align}
From the above two equations, (\ref{2.43}),  (\ref{2.431}) and (\ref{1.3-1}), it holds that
\begin{align}
\ln \frac{1-y_2(t)}{1-y_1(t)}-\ln \frac{1-y^0_2}{1-y^0_1}=&\int_0^t
\frac{b_{1}(\tau,u(\tau))}{1-y_1(\tau)}d\tau-
\int_0^t \frac{b_{2}(\tau,u(\tau))}{1-y_2(\tau)}d\tau\nonumber\\
\leq& \frac{3}{2}\ \ \mbox{for each}\ t\in [0, T_q(f^{(3)},y^0,
u)),\nonumber
\end{align}
from which we get (\ref{2.34-1}). Similarly, we can prove
(\ref{2.34-2}).

Now, we shall prove that if
\begin{align}\label{2.43-1}\lim\limits_{t\rightarrow  T_q(f^{(3)},{y}^0,u)}y_1(t)=1,\end{align}
then \begin{align}\label{2.43-2}\lim\limits_{t\rightarrow  T_q(f^{(3)},{y}^0,u)}y_2(t)=1.\end{align}

Indeed, by (\ref{2.37}) and (\ref{2.39}), it is clear that $y_2(\cdot)$
is monotonously increasing over $[0, T_q(f^{(3)},y^0,
u))$. This, combined with (\ref{2.37}), implies that
$\lim\limits_{t\rightarrow T_q(f^{(3)},y^0, u)} y_2(t)$
exists and $\lim\limits_{t\rightarrow T_q(f^{(3)},y^0, u)}
y_2(t)\leq 1$.

By contradiction, if $\lim\limits_{t\rightarrow T_q(f^{(3)},y^0, u)} y_2(t)=\beta_1<1$. Then, by the
continuity and the monotonicity of the solution  $y_2(\cdot)$, it holds that $y_2(t)\leq\beta_1$
for $t\in [0,T_q(f^{(3)},y^0, u))$. Hence,
\begin{align}
\frac{1-y_2(t)}{1-y_1(t)}\geq\frac{1-\beta_1}{1-y_1(t)},\ \ t\in [0,T_q(f^{(3)},y^0,u)),\nonumber
\end{align}
from which and (\ref{2.43-1}), it follows that
$$\lim\limits_{t\rightarrow T_q(f^{(3)},y^0, u)} \frac{1-y_2(t)}{1-y_1(t)}=+\infty,$$
which contradicts with (\ref{2.34-1}). Thus, we get (\ref{2.43-2}).

Similarly, we can prove that  if
$\lim\limits_{t\rightarrow  T_q(f^{(3)},{y}^0,u)}y_1(t)=1$,
then $\lim\limits_{t\rightarrow  T_q(f^{(3)},{y}^0,u)}y_2(t)=1$. This completes the proof of
Step 2.

{\it Step 3 To prove (\ref{2.33}) in Case 1}

We first prove that
\begin{align}\label{2.43-3}\lim\limits_{t\rightarrow T_q(f^{(3)},y^0, u)} y_1(t)=\lim\limits_{t\rightarrow T_q(f^{(3)},y^0,u)}
y_2(t)=1.\end{align}
By (\ref{2.37}), (\ref{2.38}), (\ref{2.39}), both $y_1(\cdot)$ and $y_2(\cdot)$ are monotonously increasing over
the interval $[0, T_q(f^{(3)},y^0,u))$. If (\ref{2.43-3}) did not hold, then
by Step 2 of this case and (\ref{2.37}), $\lim\limits_{t\rightarrow T_q(f^{(3)},y^0, u)} y_1(t)<1$ and
$\lim\limits_{t\rightarrow T_q(f^{(3)},y^0, u)} y_2(t)<1$. Thus by the
continuity of the solution  $y_1(\cdot)$ and $y_2(\cdot)$, there exists two numbers $\beta_2<1$ and
$\beta_3<1$ such that $y_1(t)\leq\beta_2$ and $y_2(t)\leq\beta_3$
for $t\in [0,T_q(f^{(3)},y^0, u))$, $y_1(T_q(f^{(3)},y^0, u))=\beta_2$ and $y_2(T_q(f^{(3)},y^0, u))=\beta_3$.
Then, we can extend the solution
$y(\cdot)$ and can find an interval $[T_q(f^{(3)},y^0, u), T_q(f^{(3)},y^0, u) +\delta_1)$ with
$\delta_1$ sufficiently small such that $y_1(t)<1$ and $y_2(t)<1$ on the interval
$[T_q(f^{(3)},y^0, u), T_q(f^{(3)},y^0, u)+\delta_1)$. However, $I(0,f^{(3)},y^0,u)=[0,T_q(f^{(3)},y^0,u))$ and $I(0,f^{(3)},y^0,u)$
is the maximal interval in which $y_1(t)<1$ and  $y_2(t)<1$. This is a contradiction. Thus, (\ref{2.43-3}) holds.

On the other hand, it is clear that
$$\|f^{(3)}(y(t))\|\geq \frac{1}{1-y_1(t)},\ \ t\in [0,T_q(f^{(3)},y^0,u)),$$
from which and (\ref{2.43-3}), it holds that $\lim\limits_{t\rightarrow  T_q(f^{(3)},y^0,u)}
\|f^{(3)}(y(t))\|=+\infty$. This completes the proof of (\ref{2.33}) in Case 1.

{\it Step 4. To prove (\ref{2.41-1-1}) in Case 1}

Since $y_1(t)<1$ and $y_2(t)<1$ for each $t\in [0,T_q(f^{(3)},y^0,u))$ and $y(\cdot)$ is continuous over
the interval $[0,T_q(f^{(3)},y^0,u))$, we can get $(\ref{2.41-1-1})$ in Case 1.

{\it Step 5 To prove (\ref{2.34}) in Case 1}

We shall prove (\ref{2.34}) in the following two cases.

{\it Case a. $\displaystyle\frac{1-y_1^0}{1-y_2^0}\geq 1$ in Case 1}

In this case, for each $t\in [0, T_q(f^{(3)},y^0,u))$, multiplying the first equation of
(\ref{2.42}) by $(1-y_1(\tau))^{1/2}$ and integrating it over
$[t, T_q(f^{(3)},y^0,u))$, we obtain from (\ref{2.33}) that
\begin{align}
\frac{2}{3}(1-y_1(t))^{3/2}&=\int_t^{T_q(f^{(3)},y^0,u)}
(1-y_1(\tau))^{1/2}\Big(\displaystyle\frac{1}{1-y_2(\tau)}+b_{1}(\tau,u(\tau))\Big)d\tau,
\nonumber
\end{align}
from which and (\ref{2.34-1}), we get that for each
$t\in [0, T_q(f^{(3)},{y}^0,u))$,
\begin{align}\label{2.43-4}
\frac{2}{3}(1-y_1(t))^{3/2}&\geq \int_t^{T_q(f^{(3)},y^0,u)}\Big(\displaystyle\frac{1-y^0_1}{1-y^0_2}
\frac{e^{-3/2}}{(1-y_1(\tau))^{{1/2}}}+(1-y_1(\tau))^{1/2}b_{1}(\tau,u(\tau))\Big)d\tau.
\end{align}
On the other hand, since $\displaystyle\frac{1-y^0_1}{1-y^0_2}\geq 1$,
it follows that\newpage
$$1-\frac{e^{-3/2}}{2K_0}\geq 1-\frac{e^{-3/2}(1-y^0_1)}{2K_0(1-y^0_2)},$$
from which, (\ref{2.43-4}), (\ref{2.37}) and (\ref{1.3-1}), it holds that
for each
$t\in [0, T_q(f^{(3)},{y}^0,u))$,
\begin{align}
\frac{2}{3}(1-y_1(t))^{3/2}\geq (\sqrt{2}-\frac{1}{\sqrt{2}})\sqrt{K_0}\sqrt{\frac{e^{-3/2}(1-y^0_1)}{(1-y^0_2)}}(T_q(f^{(3)},y^0,u)-t)
.\nonumber\end{align}
Using the inequality $\displaystyle\frac{1-y_1^0}{1-y_2^0}\geq 1$ again, we derive
from the above inequality that
\begin{align}
\frac{2}{3}(1-y_1(t))^{3/2}&\geq e^{-3/4}(\sqrt{2}-\frac{1}{\sqrt{2}})\sqrt{K_0}(T_q(f^{(3)},y^0,u)-t)
,\nonumber\end{align}
which implies
for each $t\in [0, T_q(f^{(3)},y^0,u))$,
\begin{align}\label{2.43-5}
\displaystyle\frac{1}{1-y_1(t)}
\leq C(T_q(f^{(3)},y^0,u)-t)^{-{2}/{3}},\end{align}
where $C$ is independent of $y^0$ and $u$.
Moreover, it follows from (\ref{2.34-2}) and (\ref{2.43-5}) that for each $t\in [0, T_q(f^{(3)},y^0,u))$,
\begin{align}
\displaystyle\frac{1}{1-y_2(t)}
&\leq \frac{e^{3/2}(1-y^0_1)}{(1-y^0_2)(1-y_1(t))}\nonumber\\
&\leq C(T_q(f^{(3)},y^0,u)-t)^{-{2}/{3}}\nonumber,\end{align}
where $C$ depends on $y^0$, but independent of $u$.
This, together with (\ref{2.43-5}), implies (\ref{2.34}) in Case 1 with $\displaystyle
\frac{1-y_1^0}{1-y_2^0}\geq 1$.

{\it Case b. $\displaystyle\frac{1-y_2^0}{1-y_1^0}>1$}

The proof of Case $b$ is very similar to that of Case $a$.
This completes the proof of (\ref{2.34}) in Case 1.

{\it Case 2. $y^0\in S^{f^{(3)}}$ with $1<y^0_1<1+\frac{e^{-3/2}}{2K_0}$, $1<y^0_2<1+\frac{e^{-3/2}}{2K_0}$}

We can use the very similar argument we used to prove Case 1 of this lemma to prove Case 2
of this lemma. This completes the proof of Lemma 2.2.3.\ \ \ \#\\
{\it {\bf Remark 2.2.3} Let $y^0\in S^{f^{(3)}}$ and $u\in \mathcal{U}_{ad}$. Since $I(0,f^{(3)},y^0,u)=[0,T_q(f^{(3)},y^0,u))$, we can conclude from the definition of $I(0,f^{(3)},y^0,u)$ that for each $t\in [0,T_q(f^{(3)},y^0,u))$,
\begin{align}
y_1(t;f^{(3)},y^0,u)<1,\ y_2(t;f^{(3)},y^0,u)<1,\ \mbox{when}\ y_1^0<1\ \mbox{and}\ y_2^0<1;\nonumber\\
y_1(t;f^{(3)},y^0,u)>1,\ y_2(t;f^{(3)},y^0,u)>1,\ \mbox{when}\ y_1^0>1\ \mbox{and}\ y_2^0>1.\nonumber
\end{align}}

\subsection{Uniform interval of non-quenching}
\ \ \ \ \ {\it Let $[0,T]$ be a time interval. Given $f\in \Lambda$, $y^0\in S^f$ and $u\in \mathcal{U}_{ad}$, we say
that the solution $y(\cdot\ ;f,y^0,u)$ does not quench on the interval $[0,T]$ if for each $t\in[0,T]$, we have
\begin{equation}
\|f(y(t;f,y^0,u))\|<+\infty.\nonumber
\end{equation}}\\
{\it {\bf Lemma 2.3}
Suppose that $f\in \Lambda$ and $[0,t_1]$ be a closed interval. Assume $y^0\in S^{f}$.
Let
$u(\cdot)$ and $\{u_k(\cdot)\}_{k=1}^\infty$ be an
element and a bounded sequence in  $\mathcal{U}_{ad}$, respectively. Suppose that
\begin{align}
u_k(\tau)\rightharpoonup u(\tau) \;\;\mbox{weakly star in}\;\;
L^{\infty}(0,t_1;\mathbb{R}^2)\ \nonumber
\end{align}
 and
the solution $y( \cdot\ ;f,
y^0,u)$ does not quench on the interval $[0, t_1]$.  Then there is a
natural number $k_0$ such that for each  $k$ with $k\geq k_0$, the
solution $y(\cdot\ ; f,y^0, u_k)$ does not quench on the interval
$[0,t_1]$.
}

{\bf Proof.}
We shall prove Lemma 2.3 in the following three cases.

{\it Case 1. $f=f^{(1)}$}

Here we only prove the case where $y^0\in S^{f^{(1)}}$ with $1-\frac{1}{2K_0}<
y^0_1<1$ and $y^0_2>K_0+\frac{1}{K_0}-1$, while we can use the very similar
 argument to prove the case where $y^0\in S^{f^{(1)}}$ with $1<
y^0_1<1+\frac{1}{2K_0}$ and $y^0_2>K_0+1$.

Let $y^0\in S^{f^{(1)}}$ with $1-\frac{1}{2K_0}<
y^0_1<1$, $y^0_2>K_0+\frac{1}{K_0}-1$.
As a matter of convenience, write $y(\cdot)$ for $y(\cdot\ ;f^{(1)},y^0,u)$, $z(\cdot)$  for ${1}/{(1-y_1(\cdot))}$
and for each $k$, write $y_k(\cdot)$ for $y(\cdot\ ;f^{(1)},y^0,u_k)$,
$z_k(\cdot)$ for ${1}/{(1-y_{1k}(\cdot))}$, respectively. Let $d_1$
be a positive number and set
$$V=\{v\in\mathbb{R}^{1}\; ;\; d(v, z[0,t_1])<d_1\}.$$
Here, $z[0, t_1]$ stands for the set $ \{ z(t) : t\in [0,
t_1]\}$ and $d(v, z[0,t_1])$ denotes the distance from the point
$v$ to the set $z[0, t_1]$ in $\mathbb{R}^{1}$.

By Lemma 2.1.1, Lemma 2.2.1 and Remark 2.2.1, both
$y(\cdot)$ and $y_k(\cdot)$ ($k=1,2\cdots$) quench at finite time,
\begin{align}\label{2.47-1}\  1-\frac{1}{2K_0}<y_1(t)<1,\ y_{2}(t)>K_0+\frac{1}{K_0}-1,\ \   t\in[0,t_1],
\end{align}
and for
each $k\ \mbox{and}\ t\in[0,T_q(f^{(1)},y^0,u_k))$,
\begin{align}\label{2.47-1-1-1}T_q(f^{(1)},y^0,u_k)<+\infty,\ & 1-\frac{1}{2K_0}<y_{1k}(t)<1, \ \ y_{2k}(t)>K_0+\frac{1}{K_0}-1,\\
& \lim\limits_{t\rightarrow T_q(f^{(1)},y^0,u_k)} y_{1k}(t)=1,\nonumber
\end{align} from which, we get
$$\lim\limits_{t\rightarrow
T_q(f^{(1)},y^0,u_k)} z_{k}(t)=+\infty.$$
Then, each solution $z_k(\cdot)$ is either in the space $C([0,t_1];
V)$ or in the space $C([0,\alpha_k); V)$, where the number
$\alpha_k$ holds the following properties:
$$(a)\ \ \ \ \ 0< \alpha_k\leq t_1;$$
$$(b)   \lim\limits_{t\rightarrow \alpha_k}z_k(t)\in
\partial V.$$

Now, we claim the following property $(B)$: {\it For all but a
finite number of $k$, solutions $z_k(\cdot)$ are in the space
$C([0,t_1]; V)$.}

 By seeking a contradiction, suppose that there
existed a subsequence of the sequence $\{z_k(\cdot)\}_{k=1}^\infty$,
still denoted in the same way, such that each solution $z_k(\cdot)$
is not in the space $C([0,t_1]; V)$. Then we would have
$z_k(\cdot)\in C([0,\alpha_k); V)$ for each  $k$. Now, corresponding
to each $k$, we define a function $e_k(\cdot)$ over $[0, \alpha_k)$
by setting
$e_k(t)=|z_k(t)-z(t)|+|y_{1k}(t)-y_1(t)|+|y_{2k}(t)-y_2(t)|$, $t\in
[0, \alpha_k)$. Then, we have
\begin{align}\label{2.48}
 e_k(t)=&\big|\frac{1}{1-y_{1k}(t)}-\frac{1}{1-y_{1}(t)}\big|
 +|y_{1k}(t)-y_1(t)|
+|y_{2k}(t)-y_2(t)|\nonumber\\
=&\frac{|y_{1k}(t)-y_1(t)|}{(1-y_{1k}(t))(1-y_1(t))}
+|y_{1k}(t)-y_1(t)|
+|y_{2k}(t)-y_2(t)|,\ \ \ \ \ \ t \in [0, \alpha_k).
\end{align}
Since for all $k$, $z_k(t)\in V$ for each $t\in [0, \alpha_k)$, and
$z(t)\in V$ for all $t\in [0, t_1]$, it holds that
\begin{align}\label{2.49}
|z(t)|\leq C,\ |z_k(t)|\leq C,\ \ t\in [0, \alpha_k).
\end{align}
Here and throughout the proof, $C$ is a positive constant independent of $k$ and $t$, which may 
be different in different context.
From (\ref{2.48}) and (\ref{2.49}), it holds that
\begin{align}
 e_k(t)\leq C\Big(|y_{1k}(t)-y_1(t)|
+|y_{2k}(t)-y_2(t)|\Big),\ \ \ \  t \in [0, \alpha_k).\nonumber
\end{align}
Then, by system (\ref{1.1}) with $f=f^{(1)}$,
it follows that for each $t\in[0,\alpha_k)$,
\begin{align}\label{2.51}  e_k(t)\leq& C\Big(|y_{1k}(t)-y_1(t)|
+|y_{2k}(\tau)-y_2(\tau)|\Big)\nonumber\\\leq &C\int_{0}^t \Big(\big|\frac{y_{2k}(\tau)}{1-y_{1k}(\tau)}-
\frac{y_{2}(\tau)}{1-y_{1}(\tau)}\big|
 +|y_{1k}(\tau)-y_1(\tau)|
+|y_{2k}(\tau)-y_2(\tau)|\Big)d\tau\nonumber\\&+C\Big\|\int_{0}^t
[B(\tau)u_k(\tau)-B(\tau)u(\tau)]d\tau\Big\|\nonumber\\\leq&C\int_0^t
\Big\{\frac{|y_{2k}(\tau)-y_2(\tau)|+|y_1(\tau)||y_{2k}(\tau)-y_2(\tau)|+|y_2(\tau)||y_{1k}(\tau)-y_1(\tau)|
}{(1-y_{1k}(\tau))(1-y_{1}(\tau))}\nonumber\\&+|y_{1k}(\tau)-y_1(\tau)|
+|y_{2k}(\tau)-y_2(\tau)|\Big\}d\tau+C\Big\|\int_{0}^t
[B(\tau)u_k(\tau)-B(\tau)u(\tau)]d\tau\Big\|.
\end{align}
On the other hand, since $y(\cdot)$ is continuous over
$[0,t_1]$, we conclude that
\begin{align}\label{2.51-1}|y_2(t)|\leq C,\ \ t\in[0,t_1].
\end{align}
Thus, by (\ref{2.47-1}), (\ref{2.49}), (\ref{2.51}) and (\ref{2.51-1}), we can  get the
following estimate:
\begin{align}e_k(t)\leq &C\int_{0}^t \Big(|y_{1k}(\tau)-y_1(\tau)|
+|y_{2k}(\tau)-y_2(\tau)|\Big)d\tau\nonumber\\&+C\Big\|\int_{0}^t
[B(\tau)u_k(\tau)-B(\tau)u(\tau)]d\tau\Big\|\nonumber,\ \ t \in [0, \alpha_k),\nonumber
\end{align}
from which and by Gronwall's inequality, it holds that
\begin{align}\label{2.52} e_k(t)\leq &\widetilde{C}\Big(\|h_k(t)\|+
\int_{0}^t \|h_k(\tau)\|d\tau\Big),\ \ t \in [0, \alpha_k).
\end{align}
Here, $\widetilde{C}$ is a positive number independent of $k$ and
$t$, and $h_k(t)=\int_{0}^t
[B(\tau)u_k(\tau)-B(\tau)u(\tau)]d\tau$,
$t \in [0, t_1]$.

On the other hand, from the properties held by $u(\cdot)$ and
$\{u_k(\cdot)\}_{k=1}^\infty$, we can easily derive that $
\lim\limits_{k\rightarrow +\infty}h_k(t)= 0$, for any $ t \in
[0,t_1]$, and that the sequence $\{h_k(\cdot)\}_{k=1}^\infty$ is
uniformly bounded and equicontinuous on the interval $[0,t_1]$.
Then, it follows from the Arzela-Ascoli theorem  that
\begin{align}\label
{2.53} h_k(\cdot)\rightarrow 0\ \;\mbox {uniformly on }\;
[0,t_1]\;\;\mbox{as}\ k\rightarrow +\infty.
\end{align}

Write $\varepsilon_0=d(\partial V, z[0,t_1])$, the distance
between the set $\partial V$ and the set $z[0,t_1]$ in
$\mathbb{R}^{1}$. From the definition of the set $V$, it follows
that $\varepsilon_0>0$. Then, we can use  (\ref{2.53}) to find  a natural
number $k^0$ enjoying the following property:
\begin{align}\label
{2.54} \widetilde{C}\Big(\| h_{k^0}(t)\|+\int_{0}^t
\|h_{k^0}(\tau)\|d\tau\Big)<\frac{\varepsilon_0}{2}\ \ \mbox{for
every}\ t \in [0, t_1].
\end{align}
Now it follows from the definition of $e_k(\cdot)$, (\ref{2.52}) and (\ref{2.54}) that
\begin{align}\label
{2.52-1}|z_{k^0}(t)-z(t)|\leq e_{k^0}(t)
<\frac{\varepsilon_0}{2}\ \ \ \mbox{for all}\ t \in [0,
\alpha_{k^0}).
\end{align}
However, according to the property $(b)$ held by $\alpha_{k^0}$, we
can find  a number $\widetilde{t}\in [0, \alpha_{k^0})$ such that
\begin{align}
d(\; z_{k^0}(\;\widetilde{t}\;), \partial
V)<\frac{\varepsilon_0}{2}.\nonumber
\end{align}
This, together with (\ref{2.52-1}), implies that
\begin{align}d( \;z(\; \widetilde{t}\; ), \partial V \;)<\varepsilon_0=d(\;\partial V,
z[0,t_1]\;),\nonumber
\end{align}
which leads to a contradiction. This completes the proof of
Property $(B)$.

By Property $(B)$ and the definition
of $V$, there exists a natural number $k_0$ such that for all $k\geq
k_0$, \begin{align}\label{2.54-1-1-2}|z_k(t)|=\frac{1}{1-y_{1k}(t)}\leq C,\ t\in[0,t_1].\end{align}
Hence, $T_q(f^{(1)},y^0,u_k)>t_1$ for each $k$. This, together with (\ref{2.47-1-1-1}), implies that
for all $k\geq
k_0$,
\begin{align}\label{2.54-1-1-1}|y_{1k}(t)|\leq C,\ t\in[0,t_1].\end{align}
On the other hand, by system (\ref{1.1}) with $f=f^{(1)}$, it holds that for all $k\geq
k_0$,
\begin{align}
y_{2k}(t)-y^0_2=\int_0^{t}(y_{1k}(\tau)+y_{2k}(\tau))d\tau+\int_0^{t}b_2(\tau,u_k(\tau))d\tau,\ t\in [0,t_1],\nonumber
\end{align}
from which, (\ref{2.54-1-1-1}) and (\ref{1.3-1}), we get that for all $k\geq
k_0$,
\begin{align}
|y_{2k}(t)|&\leq |y^0_2|+Ct_1+\int_0^{t}|y_{2k}(\tau)|d\tau,\ t\in [0,t_1].\nonumber
\end{align}
Then, by Gronwall's inequality, we derive that for all $k\geq
k_0$,
\begin{eqnarray}
|y_{2k}(t)|\leq C,\ \ \ t\in [0,t_1].\nonumber
\end{eqnarray}
This, together with (\ref{2.54-1-1-2}) and (\ref{2.54-1-1-1}) implies that for all $k\geq
k_0$,
\begin{eqnarray}
\|f^{(1)}(y_{k}(t))\|\leq C<+\infty,\ \ \ t\in [0,t_1].\nonumber
\end{eqnarray}

This completes the proof of Case 1.

{\it Case 2. $f=f^{(2)}$}

We only give a sketch proof
in the case where ${y}^0\in S^{f^{(2)}}$ with $1-\frac{1}{2K_0+1}<\|{y}^0\|<1$. The
proof in the
case where ${y}^0\in S^{f^{(2)}}$ with $1<\|{y}^0\|<1+\frac{1}{2K_0}$ is similar.
Let ${y}^0\in S^{f^{(2)}}$ with $1-\frac{1}{2K_0+1}<\|{y}^0\|<1$.
Write $y(\cdot)$ for $y(\cdot\ ;f^{(2)},y^0,u)$, $z(\cdot)$ for  $1/(1-\|y(\cdot)\|)$ and
for each $k$, 
$y_k(\cdot)$ for  $y(\cdot\ ;f^{(2)},y^0,u_k)$ and 
$z_k(\cdot)$ for  $1/(1-\|y_{k}(\cdot)\|)$, respectively. Let $d_2$
be a positive number and set
$$V=\{v\in\mathbb{R}^{1}\; ;\; d(v, z[0,t_1])<d_2\}.$$
Here, $z[0, t_1]$ stands for the set $ \{z(t) : t\in [0,
t_1]\}$ and $d(v, z[0,t_1])$ denotes the distance from the point
$v$ to the set $z[0, t_1]$ in $\mathbb{R}^{1}$. By Lemma 2.2.2,
we have $T_q(f^{(2)},y^0,u_k)<+\infty$ for each $k$. Moreover, \begin{align} \lim\limits_{t\rightarrow
T_q(f^{(2)},y^0,u_k)}\|y_k(t)\|=1\ \mbox{for each}\ k.\nonumber
\end{align}
Then,
each solution $z_k(\cdot)$ is either in the space $C([0,t_1];
V)$ or in the space $C([0,\beta_k); V)$, where the number
$\beta_k$ holds the following properties:
$$(a)\ \ \ \ \ 0<\beta_k\leq t_1;$$
$$(b)   \lim\limits_{t\rightarrow \beta_k}z_k(t)\in
\partial V.$$
In this case, we define a function $\widetilde{e}_k(\cdot)$ over $[0, \beta_k)$
by setting
$\widetilde{e}_k(t)=|z_k(t)-z(t)|+\|y_{k}(t)-y(t)\|$, $t\in
[0, \beta_k)$.
Then, we can use the similar argument we used to prove Case 1 of this lemma to complete the proof of Case 2.

{\it Case 3 $f=f^{(3)}$}

We only give a sketch proof
the case where ${y}^0\in S^{f^{(3)}}$ with $1-\frac{e^{-3/2}}{2K_0}<
y^0_1<1$ and $1-\frac{e^{-3/2}}{2K_0}<y^0_2<1$. The
proof in the
case where ${y}^0\in S^{f^{(3)}}$ with $1<{y}^0_1<1+\frac{e^{-3/2}}{2K_0}$ and
$1<{y}^0_2<1+\frac{e^{-3/2}}{2K_0}$ is similar.
Write $y(\cdot)$ and
$y_k(\cdot)$ for the solutions $y(\cdot\ ;f^{(3)},y^0,u)$
and  $y(\cdot\ ;f^{(3)},y^0,u_k)$, respectively. Write $z(\cdot)$ and
$z_k(\cdot)$ for the functions $1/(1-y_1(\cdot))+1/(1-y_2(\cdot))$
and  $1/(1-y_{1k}(\cdot))+1/(1-y_{2k}(\cdot))$, respectively.  Let $d_3$
be a positive number and set
$$V=\{v\in\mathbb{R}^{1}\; ;\; d(v, z[0,t_1])<d_3\}.$$
Here, $z[0, t_1]$ stands for the set $ \{z(t) : t\in [0,
t_1]\}$ and $d(v, z[0,t_1])$ denotes the distance from the point
$v$ to the set $z[0, t_1]$ in $\mathbb{R}^{1}$. By Lemma 2.2.3,  we have
$T_q(f^{(3)},y^0,u_k)<+\infty$ for each $k$. Moreover,
\begin{align} \lim\limits_{t\rightarrow
T_q(f^{(3)},y^0,u_k)}y_{1k}(t)=
\lim\limits_{t\rightarrow  T_q(f^{(3)},y^0,u_k)}y_{2k}(t)=1\ \mbox{for each}\ k.\nonumber
\end{align}
Then,
each solution $z_k(\cdot)$ is either in the space $C([0,t_1];
V)$ or in the space $C([0,\gamma_k); V)$, where the number
$\gamma_k$ holds the following properties:
$$(a)\ \ \ \ \ 0< \gamma_k\leq t_1;$$
$$(b)   \lim\limits_{t\rightarrow \gamma_k}z_k(t)\in
\partial V.$$
Then, we can use the similar argument we used to prove Case 1 of this lemma
to prove Case 3.\ \ \ \#

\section {Existence of Optimal Control}

\ \ \ \ \ {\bf Proof of Theorem 1.1.}
 In this section, we shall only give the proof of Theorem 1.1 in the case where $y^0\in S^{f^{(1)}}$
with $1-\frac{1}{2K_0}<y^0_1<1$ and $y^0_2>K_0+\frac{1}{K_0}-1$. We can use very similar arguments to prove Theorem 1.1 in the
cases where $y^0\in S^{f^{(1)}}$ with $1<y^0_1<1+\frac{1}{2K_0}$, $y^0_2>K_0+1$, $y^0\in S^{f^{(2)}}$
and $y^0\in S^{f^{(3)}}$.

Let $y^0\in S^{f^{(1)}}$
with $1-\frac{1}{2K_0}<y^0_1<1$, $y^0_2>K_0+\frac{1}{K_0}-1$. By Lemma 2.2.1 and Remark 2.2.1, it holds that
for each $u\in \mathcal{U}_{ad}$,
$$T_q(f^{(1)},y^0,u)<+\infty,\ \ \lim\limits_{t\rightarrow  T_q(f^{(1)},y^0,u)}y_1(t; f^{(1)},y^0, u)=1
$$and
\begin{align}\ y_1(t;f^{(1)},y^0,u)<1,\ t\in [0,T_q(f^{(1)},y^0,u)).\nonumber
\end{align}
Thus
$$
t^*=\inf\limits_{u\in \mathcal{U}_{ad}} T_q(f^{(1)},y^0,
u)<+\infty.
$$
Then, we can utilize the definitions of $t^*$ to get  a sequence
$\{u_k(\cdot)\}_{k=1}^\infty$ of  controls in the set
$\mathcal{U}_{ad}$ holding the following properties: $(1)$ Each
$t_k=:T_q(f^{(1)},y^0,u_k)$ is a positive  number; (2) $t_1\geq t_2\geq
\cdots \geq t_k\cdots\geq t^*$ and $t_k \rightarrow t^*$ as
${k\rightarrow +\infty}$; (3) $\lim\limits_{t\rightarrow t_k} y_{1}(t;f^{(1)},y^0,u_k)=1$ for all $k$.

For each $k$, write $y_k(\cdot)$ simply for $y(\cdot\ ;f^{(1)},y^0,u_k)$.
Now, we shall complete the proof by the following two steps.

{\it Step 1. To prove $t^*>0$.}

Indeed, by Lemma 2.1.1, and Remark 2.2.1, we have
\begin{align}\label{3.1}1-\frac{1}{2K_0}<y_{1k}(t)<1\ \ \mbox{for each}\ k \ \mbox{and}\ t\in [0,t_k).\end{align}
By system (\ref{1.1}) with $f=f^{(1)}$, it holds that for each $k$,
\begin{align}
y_{2k}(t)-y^0_2=\int_0^{t}(y_{1k}(\tau)+y_{2k}(\tau))d\tau+\int_0^{t}b_2(\tau,u_k(\tau))d\tau,\ t\in [0,t_k),\nonumber
\end{align}
from which, (\ref{3.1}) and (\ref{1.3-1}), we get that for each $k$,
\begin{align}
|y_{2k}(t)|&\leq |y^0_2|+Ct_k+\int_0^{t}|y_{2k}(\tau)|d\tau,\ t\in [0,t_k).\nonumber
\end{align}
Here and in what follows, $C$ is a positive constant independent of $k$ and $t$, which may be different
in different context.
Then, by the property (2) held by sequence $\{t_k\}$ and Gronwall's inequality, we derive that for each $k$,
\begin{eqnarray}\label{3.2}
|y_{2k}(t)|\leq C,\ \ \ t\in [0,t_k).
\end{eqnarray}
On the other hand, from system (\ref{1.1}) with $f=f^{(1)}$, we have that for each $k$,
\begin{align}
y_{1k}(t)-y_1^0=\int_0^{t}\frac{y_{2k}(\tau)}{1-y_{1k}(\tau)}d\tau+\int_0^{t}b_1(\tau,u_k(\tau))d\tau,\ \
t\in[0,t_k),\nonumber
\end{align}
from which, (\ref{3.2}), (\ref{2.12}) in Lemma 2.2.1 and (\ref{1.3-1}), it holds that for each $k$,
\begin{align}\label{3.3}
|y_{1k}(t)-y^0_1|&\leq \int_0^{t_k}|\frac{y_{2k}(\tau)}{1-y_{1k}(\tau)}|d\tau+\int_0^{t_k}|b_1(\tau,u_k(\tau))|d\tau\nonumber\\
&\leq C\int_0^{t_k}{(t_k-\tau)^{-{2}/{3}}}d\tau+K_0t_k,\nonumber\\&
\leq C(t_k)^{1/3}+K_0t_k, \ \ \ t\in [0,t_k).
\end{align}
If $t^*=0$, then from the property (2) held by sequence $\{t_k\}$, we have
that the right side of the above inequality tends to $0$ as $k\rightarrow+\infty$.
This, together with the inequalities (\ref{3.3}) and $y^0_1<1$, implies that we can find a natural number $K_1$ and a positive
number $\beta_0$ such that
$$y_{1K_1}(t)\leq y^0_1+\beta_0<1\ \mbox{for all}\ t\in[0,t_{K_1}),$$
which contradicts the property (3) held by sequence $\{t_{K_1}\}$. This completes the proof of Step 1.

{\it Step 2. The existence of optimal control for the problem $(P)_{y^0}^{f^{(1)}}$}

Fix such a number $T$ that $T>t^*$. It is clear that there
exist a function ${u}^*$ in $L^\infty( (0,T);\mathbb{R}^2)$ and a subsequence of the sequence
$\{u_k(\cdot)\}_{k=1}^\infty$, still denoted in the same way, such
that
\begin{align}\label{3.4} u_k\rightharpoonup {u}^*\
\mbox{weakly star in }\ L^\infty(0,T; \mathbb{R}^2) \ \ \mbox{as}
\ k\rightarrow+\infty.
\end{align}
We extend the function ${u}^*$ by setting it to be zero on
the interval $[T,+\infty)$ and denote the extension by
${u}^*$ again. Obviously,  this extended function
${u}^*$ is in the set $\mathcal{U}_{ad}$.

Now, we shall prove  that {\it $u^*$ is an optimal control for the problem $(P)_{y^0}^{f^{(1)}}$.} We shall
carry out its proof by the following two claims.

Claim One: {\it  By the definition of $t^*$, it is obvious that
the solution $y( \cdot\ ;f^{(1)},y^0,u^*)$ does not quench at any time in the interval $[0,
t^*)$.}

Claim Two:  {\it  $t^*=T_q(f^{(1)},y^0,u^*)$.} From the definition of $t^*$,
it holds that $T_q(f^{(1)},y^0,\\u^*)\geq t^*$. By seeking a contradiction, suppose that $t^*<
T_q(f^{(1)},y^0,u^*)$. Then we would find a number $\delta_0$ with $(t^*+\delta_0)<\min\{T,T_q(f^{(1)},y^0,u^*)\}$
such that the solution $y( \cdot\ ;f^{(1)},y^0,u^*)$ does not quench on the
interval $[0, t^*+\delta_0]$.
Thus, it follows from (\ref{3.4}) and  Lemma 2.3 that there exists a natural
number $\widehat{k}$ such that when $k\geq \widehat{k}$, the  solution $y( \cdot\
;f^{(1)},y^0,u_{k})$ does not quench on the interval $[0, t^*+\delta_0]$.
Thus, $t_k=T_q(f^{(1)},y^0,u_k)>t^*+\delta_0$ when $k\geq \widehat{k}$.
Now, according to (\ref{2.12}) in Lemma 2.2.1, we can easily verify
that for all $k\geq \widehat{k}$,
\begin{align}
\frac{1}{1-y_1(t;f^{(1)},y^0,u_{k})}\leq C(t_k-t)^{-2/3}\leq
(\frac{\delta_0}{2})^{-2/3},\ \ t\in[0, t^*+\frac{\delta_0}{2}],\nonumber
\end{align}
where $C$ is independent of $k$. This, together with the property
(2) held by $\{t_k\}_{k=1}^\infty$, gives a positive constant $C$
independent of $k$ such that $1/(1-y_1(t_k;f^{(1)},y^0, u_{k}))\leq
C$, for all $k\geq \widehat{k}$. This contradicts with the
property (3) by $\{t_k\}_{k=1}^\infty$.

Thus, we have complete the proof of Theorem 1.1.\ \ \ \#

\section{Proof of Pontryagin Maximum Principle}
\ \ \ \ \ \ {\bf Proof of Theorem 1.2.} We shall only give the proof of Theorem 1.2 in the case where $y^0\in S^{f^{(1)}}$
with $1-\frac{1}{2K_0}<y^0_1<1$ and $y^0_2>K_0+\frac{1}{K_0}-1$. We can use very similar arguments to prove Theorem 1.2 in the
case where $y^0\in S^{f^{(1)}}$ with $1<y^0_1<1+\frac{1}{2K_0}$ and $y^0_2>K_0+1$.

 Let $y^0\in S^{f^{(1)}}$
with $1-\frac{1}{2K_0}<y^0_1<1$ and $y^0_2>K_0+\frac{1}{K_0}-1$. We shall
prove the theorem in a series of steps as follows.

{\it Step 1. To set up a penalty functional and to study the related
properties}

Since  the number $t^*$ is the optimal time for the problem $(P)_{y^0}^{f^{(1)}}$,
the solution $y( \cdot\ ;f^{(1)},y^0, u)$, corresponding to each $u$ in the
set $\mathcal{U}_{ad}$, does not quench on $[0, t^*-\varepsilon]$,
for any $\varepsilon\in (0, t^* )$. Moreover, $T_q(f^{(1)},y^0,u) \geq t^*$
for all $u \in \mathcal{U}_{ad}$. Let $T^*$   be a fixed number such
that $T^*>t^*$. Write $\mathcal{U}[0,T^*]$ for the set $\{u|_{[0,T^*
]};\ u\in \mathcal{U}_{ad}\}$. We introduce the Ekeland distance
$d^*$ over the set $\mathcal{U}[0,T^*]$ by setting
$$d^*(u, v)=\mbox{meas} ( \{t\in [0,T^*]\ ;\ u(t)\neq v(t)\})\;\;\mbox
{for all}\; u, v\in \mathcal{U}_{ad}.$$
Here and in what follows, $\mbox{meas} (E)$ stands for the Lebesgue
measure of a measurable set $E$ in $\mathbb{R}^1$. Then
$(\mathcal{U}[0,T^*], d^*)$ forms a completed metric space (see \cite{Li-Yong},
p. 145). For each $\varepsilon\in ( 0, t^* )$, we define a penalty
functional $J_\varepsilon:(\mathcal{U}[0, T^* ], d^*)\rightarrow
R^+$ by setting
$$
J_\varepsilon(u(\cdot))=|y_1( t^*\mbox{$-$}\varepsilon; f^{(1)},y^0, u)-1|^{2}/2.
$$

We claim that {\it $J_\varepsilon$ is continuous over the space
$(\mathcal{U}[0,T^*], d^* )$}. Before moving forward to the proof
of this claim, we make the following observation, which will be
often used in what follows. Since $f^{(1)}_1(y)=y_2/(1-y_1)$,
$f_2^{(1)}(y)=y_1+y_2$,  $y=(y_1,y_2)^T\in
\mathbb{R}^2$ with $y_1\neq1$, it holds that for each $y=(y_1,y_2)^T\in
\mathbb{R}^2$ with $y_1\neq1$,\newpage
\begin{align}\label{4.1-1}
f^{(1)}_y(y)=\begin{pmatrix} \displaystyle\frac{y_2}{(1-y_1)^2}
& 1   \\
\displaystyle\frac{1}{1-y_1} & 1\end{pmatrix}.
\end{align}

Now, we come back to  prove the above claim. Let $v$ be an element and
$\{u_k(\cdot)\}_{k=1}^\infty$ be a sequence in the space
$\mathcal{U}[0,T^* ]$ such that $ d^*(u_k, v)\rightarrow 0$ as
$k\rightarrow+\infty$. Then, it is clear that
$u_k\rightarrow v$ strongly in $L^1((0, T^* );\mathbb{R}^2)$ as
$k \rightarrow +\infty$. Throughout this proof, we shall write
$y(\cdot)$ simply for $y(\cdot\ ;f^{(1)},y^0,v)$ and for each $k$, write
$y_k(\cdot)$ simply for $y(\cdot\ ;f^{(1)},y^0,u_k)$; $C$ is a constant independent
of $k$ and $t$, which may be different in different context. Since $T_q(f^{(1)},y^0,u)\geq t^*$ for all $u\in
\mathcal{U}[0,T^* ]$, it holds  from the system (\ref{1.1}) with $f=f^{(1)}$ that
\begin{align}
\label{4.1}
&|y_{1k}(t)-y_1(t)|+|y_{2k}(t)-y_2(t)|\nonumber
\\\leq&
\int_{0}^t\Big(\big|\frac{y_{2k}(\tau)}{1-y_{1k}(\tau)}-\frac{y_{2}(\tau)}{1-y_{1}(\tau)}\big|
+|y_{1k}(\tau)-y_1(\tau)|+|y_{2k}(\tau)-y_2(\tau)|\Big)d\tau\nonumber\\&+C\Big\|\int_{0}^t
[B(\tau)u_k(\tau)-B(\tau)u(\tau)]d\tau\Big\|,\ \ t \in [0,
t^*-\varepsilon].
\end{align}
From Lemma 2.1.1, Lemma 2.2.1 and Remark 2.2.1, we obtain that for each  $t\in [0,
T_q(f^{(1)},y^0,v))$,
\begin{align}
\label{4.2} 1-\frac{1}{2K_0}< y_1(t)<1,\ \ \frac{1}{1-y_1(t)}\leq
C(T_q(f^{(1)},y^0,v)-t)^{-{2}/{3}};
\end{align}
For all $k$ and for each $t\in [0,
T_q(f^{(1)},y^0,u_k))$,
\begin{align}
\label{4.3} \ 1-\frac{1}{2K_0}<y_{1k}(t)<1, \ \frac{1}{1-y_{1k}(t)}\leq
C(T_q(f^{(1)},y^0,u_k)-t)^{-{2}/{3}}.
\end{align}
By the property that $T(f^{(1)},y^0, u)\geq t^*$ for all $u\in\mathcal{U}[0,T^* ]$ again, one can
easily check that for each $k$,
$$
\max\{(T_q(f^{(1)},y^0,u_k)-t)^{-{2}/{3}},
(T_q(f^{(1)},y^0,v)-t)^{-{2}/{3}}\}\leq
\varepsilon^{-{2}/{3}},\;t\in [0, t^*-\varepsilon],$$
from which, (\ref{4.2}) and (\ref{4.3}), it holds that
\begin{align}
\label{4.2-1} \frac{1}{1-y_1(t)}\leq
C\varepsilon^{-\frac{2}{3}}\;\;\mbox{for each}\; t\in [0,
t^*-\varepsilon]
\end{align}
and for each $k$, 
\begin{align}\label{4.3-1}
\frac{1}{1-y_{1k}(t)}\leq
C\varepsilon^{-\frac{2}{3}}\;\;\mbox{for each}\; t\in [0,
t^*-\varepsilon].
\end{align}

On the other hand, because $y(\cdot)$ is continuous over $[0,t^*-\varepsilon]$, we get
\begin{align}
|y_1(t)|\leq C \ \mbox{and}\ |y_2(t)|\leq C,\ t\in
[0,t^*-\varepsilon].\nonumber
\end{align}
This, together with (\ref{4.1}), (\ref{4.2}), (\ref{4.3}), (\ref{4.2-1}) and (\ref{4.3-1}) implies that for 
each $k$,
\begin{align} &|y_{1k}(t)-y_1(t)|+|y_{2k}(t)-y_2(t)|\nonumber\end{align}\begin{align}\leq&
\int_{0}^t\Big\{\frac{|y_{2k}(\tau)-y_2(\tau)|+|y_1(\tau)||y_{2k}(\tau)-y_2(\tau)|
+|y_2(\tau)||y_{1k}(\tau)-y_1(\tau)|
}{(1-y_{1k}(\tau))(1-y_{1}(\tau))}\nonumber\\&+|y_{1k}(\tau)-y_1(\tau)|
+|y_{2k}(\tau)-y_2(\tau)|\Big\}d\tau\nonumber\\&+C\Big\|\int_{0}^t
[B(\tau)u_k(\tau)-B(\tau)u(\tau)]d\tau\Big\|\nonumber\\
\leq &C\int_{0}^t
\big(|y_{1k}(\tau)-y_1(\tau)|+|y_{2k}(\tau)
-y_2(\tau)|\big)d\tau\nonumber\\&+C\Big\|\int_{0}^{t^*-\varepsilon}
[B(\tau)u_k(\tau)-B(\tau)u(\tau)]d\tau\Big\|,\ \ t \in [0,
t^*-\varepsilon]. \nonumber\end{align}

Write $\delta(k)=
\|\int_0^{t^*-\varepsilon}B(\tau)[u_k(\tau)\mbox{$-$}v(\tau)]d\tau\|.$
Clearly, $\delta(k)\rightarrow 0$ as $k\rightarrow +\infty$.
The above inequality implies that for each $k$,
\begin{align}
&|y_{1k}(t)-y_1(t)|+|y_{2k}(t)-y_2(t)|\nonumber\\\leq
&C\int^t_0
\Big(|y_{1k}(\tau)-y_1(\tau)|+|y_{2k}(\tau)-y_2(\tau)|\Big)d\tau+
C\delta(k),\ \ t\in[0,t^*-\varepsilon]\nonumber
\end{align}
Now, we can apply Gronwall's inequality to get that as $k\rightarrow
+\infty$,
\begin{equation}\label{4.3-2}
|y_{1k}(t)-y_1(t)|+|y_{2k}(t)-y_2(t)|\rightarrow
0\;\;\mbox{uniformly in}\; t\in [0, t^*-\varepsilon].
\end{equation}
Hence, we have proved the continuity of the functional
$J_\varepsilon$.

{\it Step 2. To apply the Ekeland variational principle}

It is clear that
\begin{align}
J_\varepsilon(u^*(\cdot))=\frac{1}{2}|y_1(t^*-\varepsilon;f^{(1)},y^0,
u^*)-1|^{2}=: \sigma(\varepsilon)\rightarrow
0\;\;\mbox{as}\;\varepsilon\rightarrow 0^+, \nonumber
\end{align} and
\begin{align}
J_\varepsilon(u^*(\cdot))\leq \inf\limits_{u\in \mathcal{U}[0,T^*]}
J_\varepsilon(u(\cdot))+\sigma(\varepsilon)\;\;\mbox{for each}\;
\varepsilon\in (0,t^*). \nonumber\end{align} Then, we can utilize
Ekeland's variational principle (see, for instance, \cite{Li-Yong}, p. 136-137)
to find a control $u^\varepsilon(\cdot)\in\mathcal{U}[0,T^* ]$
enjoying the following properties:
\begin{align}
\left\{\begin{array}{ll}\label {4.4-1}d^*(u^*,
u^\varepsilon )\leq \sqrt{\sigma(\varepsilon)},\\
- \sqrt{\sigma(\varepsilon)}\ d^*(v, u^\varepsilon )\leq
J_\varepsilon(v(\cdot))-J_\varepsilon(u^\varepsilon(\cdot))\;\mbox{for
all}\; v(\cdot)\in \mathcal{U}[0,T^*].\end{array}\right.
\end{align}

Let $u(\cdot)\in \mathcal{U}[0,T^* ]$. By the variant of the
Lyapunov theorem (see, for instance, \cite{Li-Yong}, Chapter 4]), we can get,
corresponding to each  $\rho\in (0,1)$, a measurable set
$E_{\rho,\varepsilon}$ in the interval $[0,T^* ]$ such that
$\mbox{meas} ({E_{\rho,\varepsilon}} )=\rho T^*$ and
\begin{align}
\label{4.5}
\Big\|\int_{E_{\rho,\varepsilon}\cap [0,t]}
B(\tau)(u-u^\varepsilon)(\tau)d\tau-\rho\int_0^t
B(\tau)(u-u^\varepsilon)(\tau)d\tau\Big\|\leq\rho^2,\;t\in [0,T^*].
\end{align}
Now, we construct the following spike function of $u^\varepsilon$
with respect to $u$ by setting
\begin{align}\label{4.6}
u_{\rho}^\varepsilon(t)=\left\{\begin{array}{ll} u^\varepsilon(t), \
\ t\in [0,T^*]\setminus E_{\rho,\varepsilon},\\ u(t),\ \ \ t\in
E_{\rho,\varepsilon}.\end{array}\right.
\end{align}
It is obvious that the control $u_{\rho}^\varepsilon(\cdot)$ is in
$\mathcal{U}[0,T^*]$. Write $y_\rho^\varepsilon(\cdot)$ and
$y^\varepsilon(\cdot)$ simply for the solutions $y( \cdot\ ;f^{(1)},y^0,
u_\rho^\varepsilon)$ and $y( \cdot\ ;f^{(1)}, y^0,u^\varepsilon)$,
respectively. Clearly, they do not quench on the interval $[0,
t^*-{\varepsilon}]$. Set
$z_\rho^\varepsilon(t)=({y_\rho^\varepsilon(t)-
y^\varepsilon(t)})/{\rho},$ $t\in [0, t^*-\varepsilon]$. Then,
for each  $t\in [0, t^*-\varepsilon]$, it holds that
\begin{align}
\label{4.7}z_\rho^\varepsilon(t)=\int_0^t \int_0^1
\big[f_y^{(1)}(y^\varepsilon+\theta
(y_\rho^\varepsilon-y^\varepsilon))\big]^{T}(\tau)d\theta
z_\rho^\varepsilon(\tau)d\tau+\int_0^t \big\{ B(\tau)[
u_\rho^\varepsilon(\tau)-u^\varepsilon(\tau)]/{\rho}
\big\}d\tau.
\end{align}

{\it Step 3. To show the uniform convergence of the family
$\{z_\rho^\varepsilon(\cdot)\}_{\rho>0}$ on the interval $[0,
t^*-\varepsilon]$ for $\rho\rightarrow 0^+$}

It  follows from (\ref{4.6}) that $d^*(u_{\rho}^\varepsilon,
u^\varepsilon)\rightarrow 0$ as $\rho\rightarrow 0^+$. Thus, we can
use  the same argument  in the proof of (\ref{4.3-2}) to get
\begin{align}
\label{4.8}y_\rho^\varepsilon(\cdot)\rightarrow
y^\varepsilon(\cdot)\ \mbox{uniformly on}\
[0,t^*-\varepsilon]\ \mbox{as}\ \rho\rightarrow 0^+.
\end{align}
On the other hand, making use of (\ref{4.5}) and (\ref{4.6}), it holds  that for all
$t\in[0,t^*-\varepsilon]$,
\begin{align}
\label{4.9}\int_0^t \big\{
B(\tau)[u_\rho^\varepsilon(\tau)-u^\varepsilon(\tau)]/{\rho}
\big\}d\tau\mbox{$=$} \int_0^t
B(\tau)(u(\tau)-u^\varepsilon(\tau))d\tau+\|r_\rho^\varepsilon(t)\|/\rho.
\end{align}
Here, the function $r_\rho^\varepsilon(t)$ has the property:
$\|r_\rho^\varepsilon(t)\|\leq \rho^2$ for all  $t\in [0,
t^*-\varepsilon]$.

Let $z^\varepsilon(\cdot)$ be the unique solution to the following
system:
\begin{align}\label{4.10}\left\{\begin{array}{ll}&\displaystyle
\frac{dz^\varepsilon(t)}{dt}=[f_y^{(1)}(y^\varepsilon(t))]^Tz^\varepsilon(t)
+B(t)(u(t)-u^\varepsilon(t)),\ \ \ \  t\in [0, t^*-\varepsilon],\\
&z^\varepsilon(0)=0.\end{array}\right.
\end{align}
Then, by (\ref{4.7}) and (\ref{4.10}), we get the following
inequality:
\begin{align}\label{4.11}
&\|z_\rho^\varepsilon(t)-z^\varepsilon(t)\|\nonumber\\\leq&\int_0^t\Big\|
\int_0^1 [f_y^{(1)}(y^\varepsilon+\theta
(y_\rho^\varepsilon-y^\varepsilon))]^{T}(\tau)d\theta\Big\|
\big\|(z_\rho^\varepsilon(\tau)-z^\varepsilon(\tau))\big\|d\tau\nonumber\\&+\int_0^t
\Big\|\big\{\int_0^1 [f_y^{(1)}(y^\varepsilon+\theta
(y_\rho^\varepsilon-y^\varepsilon))]^{T}(\tau)d\theta-[f_y^{(1)}(y^\varepsilon)]^T(\tau)\big\}\Big\|
\big\|z^\varepsilon(\tau)\big\|d\tau\nonumber\\&+\Big\|\int_0^t
\{B(\tau)[{u_\rho^\varepsilon(\tau)-u^\varepsilon(\tau)}]/{\rho}\}d\tau-\int_0^t
B(\tau)(u(\tau)-u^\varepsilon(\tau))d\tau\Big\|, \ \ t\in [0,
t^*-\varepsilon]. \end{align} Corresponding to each $t\in [0,
t^*-\varepsilon]$, write
\begin{align}
m_\rho(t)=&\int_0^t \Big\|\big\{\int_0^1 [f_y^{(1)}(y^\varepsilon+\theta
(y_\rho^\varepsilon-y^\varepsilon))]^{T}(\tau)d\theta-[f_y^{(1)}(y^\varepsilon)]^T(\tau)\big\}\Big\|
\big\|z^\varepsilon(\tau)\big\|d\tau\nonumber\\&+\Big\|\int_0^t \big\{
B(\tau)[{u_\rho^\varepsilon(\tau)-u^\varepsilon(\tau)}]/{\rho}
\big\}d\tau-\int_0^t
B(\tau)(u(\tau)\mbox{$-$}u^\varepsilon(\tau))d\tau\Big\|.\nonumber
\end{align}
Clearly, it follows from (\ref{4.8}) and (\ref{4.9}) that $
m_\rho(\cdot)\rightarrow 0$ uniformly on $[0,t^*-\varepsilon]$ as
$\rho\rightarrow 0^+$.

Since $y^\varepsilon(\cdot)$ is continuous and  $1-\frac{1}{2K_0}<y_1^\varepsilon(t)<1$ for all
$t\in [0,t^*-\varepsilon]$ (see Lemma 2.1.1 and Remark 2.2.1), there exists a constant $\widehat{C}$ with
$0<\widehat{C}<1$ such that $y_1^\varepsilon(t)<\widehat{C}$ for each $t\in
[0,t^*-\varepsilon]$. Then, it follows from (\ref{4.8}) that when
$\rho$ is sufficiently small, we have
$y_{1\rho}^\varepsilon(t)<\widehat{C}$ for all $t\in [0,t^*-\varepsilon]$.
Hence, when
$\rho$ is  sufficiently small, it holds that for each $t\in [0,t^*-\varepsilon]$,
\begin{align}\label{4.11-1}
\frac{1}{1-\big\{y_1^\varepsilon(t)+\theta
(y_{1\rho}^\varepsilon(t)-y_1^\varepsilon(t))\big\}}\leq\frac{1}{1-\widehat{C}},\
\frac{1}{\big\{1-[y_1^\varepsilon(t)+\theta(
y_{1\rho}^\varepsilon(t)-y_1^\varepsilon(t))]\big\}^2}\leq\frac{1}{(1-\widehat{C})^2}.
\end{align}
On the other hand, from the continuity of $y^\varepsilon(\cdot)$ and
(\ref{4.8}), we get that when
$\rho$ is  sufficiently small,
$$\big|y_2^\varepsilon(t)+\theta
(y_{2\rho}^\varepsilon(t)-y_2^\varepsilon(t))\big|\leq C,\ t\in
[0,t^*-\varepsilon].$$
Here and in what follows, $C$ is a constant independent of $\rho$ and
$t$, which may be different in different context. From the above inequality, (\ref{4.1-1}) and (\ref{4.11-1}), it follows
that when $\rho$ is sufficiently small,
$$\Big\|\int_0^1 \big[f_y^{(1)}(y^\varepsilon+\theta
(y_\rho^\varepsilon-y^\varepsilon))\big]^{T}(\tau)d\theta\Big\|\leq C,\ t\in
[0,t^*-\varepsilon],$$ from which and (\ref{4.11}), it holds
from the Gronwall's inequality that when $\rho$ is sufficiently small,
$$
\|z_\rho^\varepsilon(t)-z^\varepsilon(t)\|\leq m_\rho(t)+C\int_0^t
m_\rho(\tau)d\tau, \  t\in [0, t^*-\varepsilon].$$
Hence, it
holds that $ z_\rho^\varepsilon(\cdot)\rightarrow
z^\varepsilon(\cdot)$ uniformly on $ [0,t^*-\varepsilon]$ as
$\rho\rightarrow 0^+,$
 from which,  we obtain that
\begin{align}
\label{4.12}y_\rho^\varepsilon(t)=y^\varepsilon(t)+\rho
z^\varepsilon(t)+o(\rho)\ \mbox{uniformly on}\
[0,t^*-\varepsilon].
\end{align}

{\it Step 4. To get certain necessary conditions for the control
$u^\varepsilon$}

By the second inequality of (\ref
{4.4-1}) and according to the definition of
the functional $J_\varepsilon$, we can easily verify the following
inequality:
\begin{align}
-\sqrt{\sigma(\varepsilon)}T^*\leq  [\;
\frac{1}{2}|y_{1\rho}^\varepsilon(t^*-\varepsilon)-1|^{2}-
\frac{1}{2}|y_1^\varepsilon(t^*-\varepsilon)-1|^{2}\;
]/\rho.\nonumber
\end{align}
This, together with (\ref{4.8}) and (\ref{4.12}), implies that
\begin{align}
\label{4.13}-\sqrt{\sigma(\varepsilon)}T^*&\leq
\lim\limits_{\rho\rightarrow0^+} [\;
\frac{1}{2}|y_{1\rho}^\varepsilon(t^*-\varepsilon)-1|^{2}-
\frac{1}{2}|y_1^\varepsilon(t^*-\varepsilon)-1|^{2}\;
]/\rho\nonumber\\
&=(y_1^\varepsilon(t^*-\varepsilon)-1)
z_1^\varepsilon(t^*-\varepsilon).
\end{align}
Let $\psi^\varepsilon(\cdot)$ be the unique solution for the  dual
system:\newpage
\begin{align}\label{4.18}
\left\{\begin{array}{ll}&\displaystyle\frac{d\psi^\varepsilon(t)}{dt}=-f_y^{(1)}(y^\varepsilon(t))\psi^\varepsilon(t)
,\ \ \ \  t\in [0, t^*-\varepsilon],\\
&\psi_1^\varepsilon(t^*-\varepsilon)=1-y_1^\varepsilon(t^*-\varepsilon),\ \psi_2^\varepsilon(t^*-\varepsilon)
=0.\end{array}\right.
\end{align}
Then, it follows from (\ref{4.13}) and (\ref{4.18}) that
\begin{align}
\psi_1^\varepsilon(t^*-\varepsilon)
z_1^\varepsilon(t^*-\varepsilon)\leq\sqrt{\sigma(\varepsilon)}T^*,
\nonumber\end{align} which, together with (\ref{4.10}) and
(\ref{4.18}), yields the following inequality:
\begin{align}\label{4.19}\int_0^{t^*-\varepsilon}<\psi^{\varepsilon}(\tau),
B(\tau)(u(\tau)-u^{\varepsilon}(\tau))>d\tau\leq\sqrt{\sigma(\varepsilon)}T^*.
\end{align}
We view (\ref{4.18}) and (\ref{4.19}) as  necessary conditions for
the control $u^\varepsilon$.

{\it Step 5. To obtain a uniform estimate for
$\psi^\varepsilon(\cdot)$ with $\varepsilon>0$ sufficiently small}

Since $\psi^\varepsilon(\cdot)$ solves the equation (\ref{4.18}), we
see that
\begin{align}
\psi^\varepsilon(t)=\psi^\varepsilon(t^*-\varepsilon)+\int_t^{t^*-
\varepsilon}f_y^{(1)}(y^\varepsilon(\tau))\psi^\varepsilon(\tau)d\tau
,\ \ \ \ t\in [0, t^*-\varepsilon].\nonumber
\end{align}
By Lemma 2.1.1 and Remark 2.2.1, it holds that
\begin{align}\label{4.22-2}1-\frac{1}{2K_0}<y_1^\varepsilon(t)<1,\ \ y_2^\varepsilon(t)>K_0+1/K_0-1\geq 1,\ \ t\in[0,t^*-\varepsilon].\end{align}
This, combined  with (\ref{4.1-1}), shows that for each $t\in [0,
t^*-\varepsilon]$,
\begin{align}
\|\psi^\varepsilon(t)\|&\leq(1-
y_1^\varepsilon(t^*-\varepsilon))+\int_t^{t^*-
\varepsilon}\|f_y^{(1)}(y^\varepsilon(\tau))\|\ \|\psi^\varepsilon(\tau)\|d\tau\nonumber\\
&\leq(1- y_1^\varepsilon(t^*-\varepsilon))+\int_t^{t^*-
\varepsilon}\big(\frac{y_2^\varepsilon(\tau)}{(1-y_1^\varepsilon(\tau))^2}+\frac{1}{1-y_1^\varepsilon
(\tau)} +2\big)\|\psi^\varepsilon(\tau)\|d\tau.\nonumber
\end{align}
Now we can apply  the Gronwall's inequality to get that for all
$t\in [0, t^*-\varepsilon]$,
\begin{align}
\label{4.16}\|\psi^\varepsilon(t)\|
\leq(1-
y_1^\varepsilon(t^*-\varepsilon))\exp\Big\{
\int_t^{t^*-\varepsilon}\big(\frac{y_2^\varepsilon(\tau)}{(1-y_1^\varepsilon
(\tau))^2}+\frac{1}{1-y_1^\varepsilon(\tau)}
+2\big)d\tau\Big\}.
\end{align}

On the other hand, according to the system
satisfied by $y^\varepsilon(\cdot)$, it follows that for
every $t\in [0, t^*-\varepsilon]$,
\begin{align}
&\int_t^{t^*-\varepsilon}\frac{d(1-y_1^\varepsilon(\tau))}{1-y_1^\varepsilon(\tau)}\nonumber\\
=&\
\ln\frac{1-y_1^\varepsilon(t^*-\varepsilon)}{1-y^\varepsilon_1(t)}=\int_t^{t^*-\varepsilon}
-\frac{1}{1-y_1^\varepsilon(\tau)}\big(\frac{y_2^\varepsilon(\tau)}{1-y_1^\varepsilon(\tau)}+b_1(\tau,
u^\varepsilon(\tau)\big)d\tau,\nonumber
\end{align}
namely, we have the equation as follows:
\begin{align}
&1-y_1^\varepsilon(t^*-\varepsilon)\nonumber\\=&(1-y_1^\varepsilon(t))\cdot
\exp \Big\{\int_t^{t^*-\varepsilon}
-\frac{y_2^\varepsilon(\tau)}{(1-y_1^\varepsilon(\tau))^2}d\tau\Big\}\nonumber\\&\cdot
\exp\Big\{\int_t^{t^*-\varepsilon} \frac{-b_1(\tau,u^\varepsilon(\tau))}{1-y_1^\varepsilon(\tau)}d\tau\Big\}, \ \
t\in [0, t^*-\varepsilon]. \nonumber\end{align} Then,
making use of the above equation and by (\ref{4.16}), we obtain that
\begin{align}\label{4.17}\|\psi^\varepsilon(t)\|\leq(1-
y_1^\varepsilon(t))\exp\Big\{
\int_t^{t^*-\varepsilon}\big(\frac{1}{1-y_1^\varepsilon(\tau)}
+2-\frac{b_1(\tau,u^\varepsilon(\tau))}{1-y_1^\varepsilon(\tau)}\big)d\tau\Big\}, \  \
t\in [0, t^*-\varepsilon].\end{align} On the other hand,
by Lemma 2.2.1, we obtain that
$$
\frac{1}{1-y_1^\varepsilon(t)}\leq  C(T_q(f^{(1)},y^0,u^\varepsilon)-t)^{-{2}/{3}}\;\;\mbox{for each}\; t\in [0, T_q(f^{(1)},y^0,u^\varepsilon)),
$$
where $C$ is independent of $\varepsilon$ and $t$,
from which and the inequality $T_q(f^{(1)},y^0,u^\varepsilon)\geq t^*$, it follows that
\begin{align}\label{4.22-1}
\frac{1}{1-y_1^\varepsilon(t)}\leq  C(t^*-t)^{-{2}/{3}}\;\;\mbox{for each}\; t\in [0, t^*-
\varepsilon].
\end{align}
This, together with (\ref{4.17}) and (\ref{1.3-1}), implies that
\begin{align}\label{4.22}\|\psi^\varepsilon(t)\|\leq C_1(1-
y_1^\varepsilon(t)), \  \
t\in [0, t^*-\varepsilon],
\end{align}
where  $C_1$ is independent of $\varepsilon\in (0,t^*)$ and $t$.

{\it Step 6. Convergence of a subsequence of the family
$\{\psi^\varepsilon(\cdot)\}_{\varepsilon>0}$ }

First of all, corresponding to each $\varepsilon\in (0, t^*)$, we
extend the function $\psi^\varepsilon(\cdot)$ by setting $\psi_1^\varepsilon(\cdot)$ to be
$1-y_1^\varepsilon( t^*-\varepsilon )$ on $(t^*-\varepsilon,t^* ]$, and setting
$\psi_2^\varepsilon(\cdot)$ to be
0 on $(t^*-\varepsilon,t^* ]$,
denote the extended function by $\psi^\varepsilon(\cdot)$ again.
Clearly, this extended function is continuous on $[0, t^* ]$.

Now, we take a sequence  $\{\delta_m\}_{m=1}^\infty$ of numbers from
the interval $ (0,t^* )$  such that (i)
$\lim\limits_{m\rightarrow+\infty}\delta_m= 0$; (ii)
$\delta_1>\delta_2>\cdots $.  Corresponding to the number
$\delta_1$, we can take a sequence $\{\varepsilon_n\}_{n=1}^\infty$
from the set $\{\varepsilon\}_{0<\varepsilon<t^*}\;$ such that
$\lim\limits_{n\rightarrow +\infty}\varepsilon_n= 0$  and $[0,
t^*-\delta_1]\subset[0, t^*-\varepsilon_n]$ for all
$n=1,2,\cdots$.

By (\ref{4.22-2}), the sequence
$\{1-y_1^{\varepsilon_n}(\cdot)\}_{n=1}^\infty$ is uniformly bounded on
the interval $[0, t^*-\delta_1]$. This, together with (\ref{4.22}), implies that $\{\psi^{\varepsilon_n}(\cdot)\}_{n=1}^\infty$
 is uniformly bounded on $[0, t^*-\delta_1]$.

On the other hand, by system
(\ref{1.1}) with $f=f^{(1)}$, we can use the similar argument we used to prove (\ref{3.2}) to conclude that
\begin{align}\label{4.22-3}|y_2^{\varepsilon_n}(t)|\leq C,\ \ t\in[0,t^*-\delta_1],
\end{align}
where $C$ is independent of $n$ and $t$.

Next, we shall prove that the sequence
$\{\psi^{\varepsilon_n}(\cdot)\}_{n=1}^\infty$ is  equicontinuous on $[0, t^*-\delta_1]$.

Indeed, for each $s_1$, $s_2$ in the interval $[0,t^*-\delta_1]$ with $s_1<s_2$,  it follows from
(\ref{4.18}), (\ref{4.1-1}) and (\ref{4.22}) that
\begin{align}&\big|\psi^{\varepsilon_n}(s_1)-\psi^{\varepsilon_n}(s_2)\big|\nonumber\\\leq&
\int_{s_1}^{s_2}\big\|f_y^{(1)}(y^{\varepsilon_n}(t))\big\|\big\|\psi^{\varepsilon_n}(t)\big\|dt\nonumber\\\leq&
C_1\int_{s_1}^{s_2}\Big(\frac{y_2^{\varepsilon_n}(t)}{(1-y_1^{\varepsilon_n}(t))^2}+
\frac{1}{1-y_1^{\varepsilon_n}(t)}+2\Big)(1-y_1^{\varepsilon_n}(t))dt,\nonumber
\end{align}
from which, (\ref{4.22-2}), (\ref{4.22-1}) and (\ref{4.22-3}), it holds that
 \begin{align}&\big|\psi^{\varepsilon_n}(s_1)-\psi^{\varepsilon_n}(s_2)\big|\nonumber\\\leq&
C\int_{s_1}^{s_2}
\Big(\frac{1}{1-y_1^{\varepsilon_n}(t)}+1\Big)dt,\nonumber\\\leq&
C\int_{s_1}^{s_2}(t^*-t)^{-2/3}dt+C(s_2-s_1),\nonumber\\\leq&
C\big[(t^*-s_1)^{1/3}-(t^*-s_2)^{1/3}\big]+C(s_2-s_1),\nonumber
\end{align}
where $C$ is independent of $n$, $s_1$ and $s_2$. This implies that
the sequence $\{\psi^{\varepsilon_n}(\cdot)\}_{n=1}^\infty$ is
equicontinuous on $[0, t^*-\delta_1]$. Hence, we can utilize
the Arzela-Ascoli theorem to take  a subsequence
$\{\psi^{\varepsilon_{n,1}}(\cdot)\}_{n=1}^\infty$ from the sequence
$\{\psi^{\varepsilon_n}(\cdot)\}_{n=1}^\infty$ such that it is
uniformly convergent on the interval $[0, t^*-\delta_1]$.
Following this process,  we can take, corresponding to each number
$\delta_m$, a subsequence
$\{\psi^{\varepsilon_{n,m}}(\cdot)\}_{n=1}^\infty$ holding the
properties as follows: $(a_1)$
$\{\psi^{\varepsilon_{n,j}}(\cdot)\}_{n=1}^\infty\subset
\{\psi^{\varepsilon_{n,j-1}}(\cdot)\}_{n=1}^\infty$ for all
$j=2,\cdots,m$; $(a_2)$ The sequence
$\{\psi^{\varepsilon_{n,m}}(\cdot)\}_{n=1}^\infty$ is uniformly
convergent on $[0, t^*-\delta_m]$. Now, by the standard diagonal
argument, we see that the subsequence
$\{\psi^{\varepsilon_{n,n}}(\cdot)\}$ is uniformly convergent on
$[0, t^*-\delta]$ for each  $\delta\in (0, t^* )$.

Let
$\psi(t)=\lim\limits_{n\rightarrow+\infty}\psi^{\varepsilon_{n,n}}(t)$,
$t\in[0, t^*)$. Then it holds that
\begin{align}\label{4.23}\psi^{\varepsilon_{n,n}}(\cdot)\rightarrow
\psi(\cdot)\ \mbox{uniformly  on}\ [0, t^*-\delta]\ \
\mbox{as}\ n\rightarrow +\infty,\;\;\mbox{for each}\; \delta\in (0,
t^* ) .
\end{align}

{\it Step 7. To extend the function $\psi$  over the interval $[0,
t^*]$}

Since all solutions $y^{\varepsilon_{n,n}}(\cdot)$,
$n=1,2,\cdots$, does not quench on every interval $[0,
t^*-{\delta}]$ with $\delta\in(0,t^*)$ and
$\lim\limits_{n\rightarrow+\infty}d^*(u^{\varepsilon_{n,n}}, u^*)=
0$,
we can use the same arguments in the proof of (\ref{4.3-2}) to get that
\begin{align}\label{4.25}y^{\varepsilon_{n,n}}(\cdot)\rightarrow
y^*(\cdot)  \  \mbox{uniformly on}\ [0, t^*-\delta] \ \
\mbox{as}\ n\rightarrow +\infty,\;\;\mbox{for each}\; \delta\in (0,
t^* ).
\end{align}
This, together with (\ref{4.22}) and (\ref{4.23}),  yield the
following inequality:
\begin{align}\|\psi({t})\|
\leq C_1 (1-y_1^*({t}))\;\; \mbox{for all}\; t\in [0, t^* ).
\nonumber\end{align} Since  the constant $C_1$ is independent of
${t}$, it holds that
\begin{align}\label{4.25-1}
0\leq \lim\limits_{t\rightarrow t^*}\|\psi(t)\|\leq
\lim\limits_{t\rightarrow t^*}C_1 (1-y_1^*({t})).
\end{align}
On the other hand, because $t^*=T_q(f^{(1)},y^0,u^*)$,
we conclude from Lemma 2.2.1 that $\lim\limits_{t\rightarrow t^*}y_1^*({t})=1$,
from which and (\ref{4.25-1}), we obtain that $
\lim\limits_{t\rightarrow{t^*}}\psi(t)=0.$ Now, we extend the
function $\psi(\cdot)$ by setting it to be zero at the time $t^*$
and still denote the extension by $\psi(\cdot)$. Clearly, this
extended function $\psi(\cdot)$ is continuous on the interval $[0,
t^*]$ and has the property: $\psi(t^*)=0$.

{\it Step 8. To verify that the function  $\psi(\cdot)$ solves the
system (\ref{1.3})}

Clearly, the function $\psi^{\varepsilon_{n,n}}(\cdot)$ holds the
following property:
\begin{align}\label{4.26}
\psi^{\varepsilon_{n,n}}(t)=\psi^{\varepsilon_{n,n}}(t^*-\varepsilon_{n,n})+\int_t^{t^*-
\varepsilon_{n,n}}f_y^{(1)}(y^{\varepsilon_{n,n}}(\tau))\psi^{\varepsilon_{n,n}}(\tau)d\tau ,\ \
\ \ t\in [0, t^*-\varepsilon_{n,n}],
\end{align}
where $\psi^{\varepsilon_{n,n}}_1(t^*-\varepsilon_{n,n})=1-y_1^{\varepsilon_{n,n}}(t^*-\varepsilon_{n,n})$,
$\psi_2^{\varepsilon_{n,n}}(t^*-\varepsilon_{n,n})=0$.
We first claim that
\begin{align}\label{4.27}\frac{1}{1-y_1^{\varepsilon_{n,n}}(t^*
-{{\varepsilon_{n,n}}})}\rightarrow +\infty\; \  \mbox{as}
\; n\rightarrow+\infty,
\end{align}
or equivalently,
\begin{align}\label{4.27-1-1-1}\psi^{\varepsilon_{n,n}}(t^*-\varepsilon_{n,n})\rightarrow 0\; \  \mbox{as}
\; n\rightarrow+\infty.
\end{align}
If (\ref{4.27}) were not true, then there would exist a positive constant
$\beta$ such that
\begin{equation}\label{4.28}
\frac{1}{1-y_1^{\varepsilon_{n,n}}(t^*
-{{\varepsilon_{n,n}}})}\leq
\beta\;\;\mbox{for infinitely many} \; n.
\end{equation}
On the other hand, we
take a number $\gamma_1$ with
$1-\frac{1}{2K_0}<\gamma_1<1$ and $1/(1-\gamma_1)>\beta$. Because $\lim\limits_{t\rightarrow t^*}y_1^*({t})=1$,
we can find a number $\delta \in (0, t^*)$
such that $y_1^*(t^* -\delta)>\gamma_1$.
Since
$\lim\limits_{n\rightarrow+\infty}\varepsilon_{n,n}=0$, there is a natural number $N_1$ such that
$$
(t^*-\varepsilon_{n,n})>(t^*-\delta)\;\;\mbox{for
all}\;\; n\geq N_1.
$$
Now, we can utilize Lemma 2.1.1 to get the following inequality:
$$y_1^{\varepsilon_{n,n}}(t^* - \varepsilon_{n,n})\geq \gamma_1\;\;\mbox{for
all}\;\; n\geq N_1,$$
from which, we have
$$
\frac{1}{1-y_1^{\varepsilon_{n,n}}(t^* - \varepsilon_{n,n})}\geq
\frac{1}{1-\gamma_1}>\beta\;\;\mbox{for all}\; n\geq N_1,
$$
which contradicts to (\ref{4.28}). Therefore, we have proved
(\ref{4.27}).

Next, we claim that
\begin{align}\label{4.29}\lim\limits_{n\rightarrow +\infty}\int_t^{t^*-{\varepsilon_{n,n}}}
f_y^{(1)}(y^{\varepsilon_{n,n}}(\tau))\psi^{\varepsilon_{n,n}}(\tau)d\tau=
\int_t^{t^*} f_y^{(1)}(y^*(\tau))\psi(\tau)d\tau, \ \ t\in [0,t^* ).
\end{align}
Corresponding to  each $n$, we define a function
$F_n(\cdot)$ by setting
$$
F_n(t)=\left\{\begin{array}{ll}
f_y^{(1)}(y^{\varepsilon_{n,n}}(t))\psi^{\varepsilon_{n,n}}(t), & t\in [0,
t^*-\varepsilon_{n,n}],\\
0, & t\in (t^*-\varepsilon_{n,n}, t^*].
\end{array}\right.
$$
It is clear that all  functions $F_n(\cdot)$, $n=1,2,\cdots$, are
measurable on the interval $[0,t^*]$.  We shall first give an estimate on the
sequence $\{\|F_n(\cdot)\|\}_{n=1}^\infty$.  Let  $t\in [0, t^* ]$. In the case that
$n$ is such that $t\in [0, t^*-\varepsilon_{n,n}]$, by (\ref{4.1-1}), (\ref{4.22-2}), (\ref{4.22-1}), (\ref{4.22})
 and (\ref{4.22-3}), we have the following
estimate:
\begin{align}\label{4.30} \|F_n(t)\|
&=\big\|f_y^{(1)}(y^{\varepsilon_{n,n}}(t))\psi^{\varepsilon_{n,n}}(t)\big\|\nonumber\\&\leq
C\Big(\frac{|y_2^{\varepsilon_{n,n}}(t)|}{(1-y_1^{\varepsilon_{n,n}}(t))^2}+
\frac{1}{1-y_1^{\varepsilon_{n,n}}(t)}
+2\Big)(1-y_1^{\varepsilon_{n,n}}(t))
\nonumber\\&\leq C+C(t^*-t)^{-2/3}, \ t\in [0,t^*-\varepsilon_{n,n}],
\end{align}
where $C$ is
independent of $n$ and $t$.

On the other hand, if $n$ is such that
$t\in (t^*-\varepsilon_{n,n}, t^*]$, then $F_n(t)=0$. This,
together with (\ref{4.30}) implies  that
\begin{align}\label{4.30-1}
\|F_n(t)\|\leq C+C(t^*-t)^{-2/3},\;\; \mbox{for all}\; t\in [0,
t^* ]\;\mbox{ and for all }\; n=1,2,\cdots.
\end{align}
Then, we are going to show that the sequence
$\{F_n(t)\}_{n=1}^\infty$ is convergent for each $t\in [0, t^*)$.
Indeed,  corresponding to each  $t\in [0, t^* )$, there exists a
natural number $N_2$ such that $t\in  [0,
t^*-\varepsilon_{n,n}]$ when $n\geq N_2$. Then, by (\ref{4.23}) and
(\ref{4.25}), we obtain that
\begin{align}F_n(t)\rightarrow
f_y^{(1)}(y^*(t) )\psi (t)\;\; \mbox{as}\; n\rightarrow +\infty,\nonumber
\end{align}
from which and (\ref{4.30-1}), we can apply Lebesgue dominated convergence theorem to get
that
$$
\lim\limits_{n\rightarrow+\infty}\int^{t^*}_t F_n(\tau) d\tau =
\int^{t^*}_t f_y^{(1)}(y^*(\tau))\psi (\tau)d\tau,\;\;\mbox{for each}\;
t\in [0, t^* ),
$$
from which, (\ref{4.29}) follows immediately.

Now, let $s\in [0, t^*  )$. Clearly, it holds that $s\in [0,
t^*-\varepsilon_{n,n})$ for $n$ sufficiently large. Thus,
making use of (\ref{4.27-1-1-1}) and (\ref{4.29}), we can pass to the limit for
$n\rightarrow+\infty$ in (\ref{4.26}), where $t=s$, to get that
$$
\psi (s)= \int_s^{t^*} f_y^{(1)}(y^*(\tau))\psi(\tau)d\tau.
$$
Since $s$ can be  arbitrarily taken from $[0, t^* )$ in the above
equation, we already hold the first  equation of (\ref{1.3}). On the
other hand, we proved that $\psi (t^*)=0$ in the end of Step 7.
Hence, the function $\psi(\cdot)$ solves the system (\ref{1.3}).

{\it Step 9. To prove the equation (\ref{1.4})}

By the first inequality of (\ref{4.4-1}), we can take  a subsequence from
the sequence $\{u^{\varepsilon_{n,n}}(\cdot)\}_{n=1}^\infty$, still
denoted in the same way, such that when $n\rightarrow+\infty$,
\begin{equation}\label{4.32}
u^{\varepsilon_{n,n}}(t)\rightarrow u^*(t),\;\; \mbox{for almost
every }\; t\in [0, T^* ].
\end{equation}
Corresponding to  each $n$, we define a function
$H_n(\cdot)$ by setting
$$
H_n(t)=\left\{\begin{array}{ll} \psi^{\varepsilon_{n,n}}(t), & t\in
[0,
t^*-\varepsilon_{n,n}\;],\\
0, & t\in (t^*-\varepsilon_{n,n}, t^*\; ].
\end{array}\right.
$$
Then, by (\ref{4.23}) and (\ref{4.32}), it holds that for each $u(\cdot)$ in the
space $\mathcal{U}[0,T^* ]$,
$$
\lim\limits_{n\rightarrow+\infty}<H_n(t) , B(t)(
u(t)-u^{\varepsilon_{n,n}}(t))> = <\psi(t),
B(t)(u(t)-u^*(t))>\;\;\mbox{for a.e.}\;t\in [0, t^* ].
$$
Moreover, making use of (\ref{4.22-2}) and (\ref{4.22}), we
obtain that
$$
|< H_n(\tau) , B(\tau) ( u(\tau)-
u^{\varepsilon_{n,n}}(\tau))
>|\leq C\;\mbox{for a.e.}\;\tau \in [0, t^*],
$$
where $C$ is independent of $n$.
Now, we can apply Lebesgue dominated convergence theorem to get that
$$
\lim\limits_{n\rightarrow+\infty}\int_0^{t^*} < H_n(\tau) , B(\tau)
( u(\tau)-u^{\varepsilon_{n,n}}(\tau) ) > d\tau= \int_0^{t^*} <\psi
(\tau), B(\tau)( u(\tau) -u^*(\tau) )>d\tau.
$$
Namely,
\begin{align}
&\lim\limits_{n\rightarrow+\infty}\int_0^{t^*-\varepsilon_{n,n}} <
\psi^{\varepsilon_{n,n}}(\tau) , B(\tau) (
u(\tau)-u^{\varepsilon_{n,n}}(\tau) ) > d\tau \nonumber\\
 &=  \int_0^{t^*} <\psi (\tau), B(\tau)( u(\tau)
-u^*(\tau) )>d\tau,\nonumber
\end{align}
which, together with (\ref{4.19}) where
$\varepsilon=\varepsilon_{n,n}$, implies the following inequality:
$$
\int_0^{t^*} <\psi (\tau), B(\tau)( u(\tau) -u^*(\tau)
>d\tau\leq 0.
$$
Since the above inequality holds for every control $u(\cdot)\in
\mathcal{U}[0,T^* ]$, we can use the standard argument
(see, for instance, \cite{Li-Yong}, p. 157-158) to derive the equation (\ref{1.4}).

{\it Step 10. To show the non-triviality of the function
$\psi(\cdot)$}

By the equations satisfied by $y^{\varepsilon_{n,n}}(\cdot)$ and
$\psi^{\varepsilon_{n,n}}(\cdot)$, respectively,  we see that for each $t\in [0, t^*-\varepsilon_{n,n}]$,
\begin{align}
\label{4.33}
&\int_t^{t^*-\varepsilon_{n,n}}\frac{d}{d\tau}\Big(\frac{\psi_1^{\varepsilon_{n,n}}
(\tau)}{1-y_1^{\varepsilon_{n,n}}(\tau)}\Big)
d\tau\nonumber\\=
&\int_t^{t^*-\varepsilon_{n,n}}\frac{1}{1-y_1^{\varepsilon_{n,n}}(\tau)}
\Big\{-\frac{y_2^{\varepsilon_{n,n}}(\tau)}{(1-y_1^{\varepsilon_{n,n}}(\tau))^2}
\psi_1^{\varepsilon_{n,n}}(\tau)-
{\psi_2^{\varepsilon_{n,n}}(\tau)}\Big\}d\tau
\nonumber\\&+
\int_t^{t^*-\varepsilon_{n,n}}\frac{1}{(1-y_1^{\varepsilon_{n,n}}(\tau))^2}
\Big\{\frac{y_2^{\varepsilon_{n,n}}(\tau)}{1-y_1^{\varepsilon_{n,n}}(\tau)}
+b_1(\tau,u^{\varepsilon_{n,n}}(\tau))\Big\}\psi_1^{\varepsilon_{n,n}}(\tau)d\tau\nonumber\\
=&\int_t^{t^*-\varepsilon_{n,n}}\frac{-
\psi_2^{\varepsilon_{n,n}}(\tau)}{1-y_1^{\varepsilon_{n,n}}(\tau)}d\tau+
\int_t^{t^*-\varepsilon_{n,n}}\frac{\psi_1^{\varepsilon_{n,n}}(\tau)}
{(1-y_1^{\varepsilon_{n,n}}(\tau))^2}
b_1(\tau,u^{\varepsilon_{n,n}}(\tau))d\tau.
\end{align}
Because
\begin{align}\frac{\psi_1^{\varepsilon_{n,n}}(t^*-\varepsilon_{n,n})}{1-y_1^{\varepsilon_{n,n}
}(t-\varepsilon_{n,n})}=1,\nonumber
\end{align}
it follows from (\ref{4.33}) that for each $\ t\in [0, t^*-\varepsilon_{n,n}]$,
\begin{align}\label{4.34}
&1-\frac{\psi_1^{\varepsilon_{n,n}}(t)}{1-y_1^{\varepsilon_{n,n}}(t)}\nonumber\\=
&\int_t^{t^*-\varepsilon_{n,n}}\frac{-
\psi_2^{\varepsilon_{n,n}}(\tau)}{1-y_1^{\varepsilon_{n,n}}(\tau)}d\tau+
\int_t^{t^*-\varepsilon_{n,n}}\frac{\psi_1^{\varepsilon_{n,n}}(\tau)}
{(1-y_1^{\varepsilon_{n,n}}(\tau))^2}b_1(\tau,u^{\varepsilon_{n,n}}(\tau))d\tau.
\end{align}
By (\ref{4.22}), we can make use of the very similar  same arguments as those in the
proof of (\ref{4.29}) to verify that
\begin{align}
&\lim\limits_{n\rightarrow+\infty}\Big(\int_t^{t^*-\varepsilon_{n,n}}\frac{-
\psi_2^{\varepsilon_{n,n}}(\tau)}{1-y_1^{\varepsilon_{n,n}}(\tau)}d\tau+
\int_t^{t^*-\varepsilon_{n,n}}\frac{\psi_1^{\varepsilon_{n,n}}(\tau)}
{(1-y_1^{\varepsilon_{n,n}}(\tau))^2}
b_1(\tau,u^{\varepsilon_{n,n}}(\tau))d\tau\Big)\nonumber\\&=
\int_t^{t^*}\frac{-
\psi_2(\tau)}{1-y_1^*(\tau)}d\tau+
\int_t^{t^*}\frac{\psi_1(\tau)}
{(1-y_1^{*}(\tau))^2}
b_1(\tau,u^*(\tau))d\tau,\;\;\mbox{for each}\; t\in [0, t^*
),\nonumber
\end{align}
which, together (\ref{4.23}), (\ref{4.25}) and (\ref{4.34}), implies that for each
$t\in [0, t^*)$,
\begin{align}
&1-\frac{\psi_1{(t)}}{1-y_1^*(t)}= \int_t^{t^*}\frac{-
\psi_2(\tau)}{1-y_1^*(\tau)}d\tau+
\int_t^{t^*}\frac{\psi_1(\tau)} {(1-y_1^{*}(\tau))^2}
b_1(\tau,u^*(\tau))d\tau.
\nonumber
\end{align}
This shows that the function $\psi(\cdot)$ is not trivial.

{\it Step 11. To prove (\ref{1.4-1}})

Indeed, by (\ref{2.13-1}) in Lemma 2.2.1, (\ref{1.4-1}) holds.

Thus, we complete the proof of Theorem 1.2 in the case where $y^0\in S^{f^{(1)}}$
with $1-\frac{1}{2K_0}<y^0_1<1$ and $y^0_2>K_0+\frac{1}{K_0}-1$.
\ \ \ \#

{\bf Sketch proof of Theorem 1.3.} We can
use the similar argument of the proof of Theorem 1.2 to prove
Theorem 1.3.
We shall only give the key steps of the proof of this theorem in the case where $y^0\in S^{f^{(2)}}$  with $1-\frac{1}{2K_0+1}<\|y^0\|<1$.
The proof of the case where $y^0\in S^{f^{(2)}}$ with $1<\|{y}^0\|<1+\frac{1}{2K_0}$
is similar.

Since $f^{(2)}(y)=y/(1-\|y\|)$, $y=(y_1,y_2)^T\in
\mathbb{R}^2$ with $\|y\|\neq1$,
it holds that for each $y=(y_1,y_2)^T\in
\mathbb{R}^2$ with $\|y\|\neq1$,
\begin{align}\label{4.37}
f^{(2)}_y(y)=\frac{1}{(1-\|y\|)^2}\begin{pmatrix} 1-\|y\|+\displaystyle\frac{y_1^2}{\|y\|}
& \displaystyle\frac{y_1y_2}{\|y\|}   \\
\displaystyle\frac{y_1y_2}{\|y\|}& 1-\|y\|+\displaystyle\frac{y_2^2}{\|y\|}\end{pmatrix}.
\end{align}

Let $T^*>t^*$. For each $\varepsilon\in(0,t^*)$, we define a penalty
functional $J_\varepsilon:(\mathcal{U}[0, T^* ], d^*)\rightarrow
R^+$ (See p. 28 for the definition of $\mathcal{U}[0, T^* ]$ and $d^*$)  by setting
$$
J_\varepsilon(u(\cdot))=(\|y( t^*-\varepsilon;f^{(2)},y^0,u)\|-1)^2/2.
$$
Write $y^\varepsilon(\cdot)
$ for $y(\cdot\;;f^{(2)},y^0,u^\varepsilon)$. We can use the similar argument we used of the proof of
Theorem 1.2 to obtain that
$J_\varepsilon$ is  continuous over the space
$(\mathcal{U}[0,T^*], d^* )$, and to get that there exists a control $u^\varepsilon(\cdot)$
enjoys the property (\ref{4.4-1}) and the following inequality,
\begin{align}\int_0^{t^*-\varepsilon}<\psi^{\varepsilon}(\tau),
B(\tau)(u(\tau)-u^{\varepsilon}(\tau))>d\tau\leq\sqrt{\sigma(\varepsilon)}T^*,\nonumber
\end{align}
where
\begin{align}\label{4.38}
\left\{\begin{array}{ll}&\displaystyle\frac{d\psi^\varepsilon(t)}{dt}=-f_y^{(2)}(y^\varepsilon(t))\psi^\varepsilon(t)
,\ \ \ \  t\in [0, t^*-\varepsilon],\\
&\psi^\varepsilon(t^*-\varepsilon)=\displaystyle\frac{1-\|y^\varepsilon(t^*-\varepsilon)\|}
{\|y^\varepsilon(t^*-\varepsilon)\|}y^\varepsilon(t^*-\varepsilon).\end{array}\right.
\end{align}

Now, we shall obtain a uniform estimate for $\psi^\varepsilon(\cdot)$ with
$\varepsilon>0$ sufficiently small in this case.

Indeed, since $\psi^\varepsilon(\cdot)$ solves the equation
(\ref{4.38}), we see that
\begin{align}
\psi^\varepsilon(t)=\displaystyle\frac{1-\|y^\varepsilon(t^*-\varepsilon)\|}
{\|y^\varepsilon(t^*-\varepsilon)\|}y^\varepsilon(t^*-\varepsilon)+\int_t^{t^*-
\varepsilon}f_y^{(2)}(y^\varepsilon(\tau))\psi^\varepsilon(\tau)d\tau
,\ \ \ \ t\in [0, t^*-\varepsilon].\nonumber
\end{align}
By Lemma 2.1.2 and Remark 2.2.2, it holds that
\begin{align}\label{4.38-1}1-\frac{1}{2K_0+1}<\|y^\varepsilon(t)\|<1,\ \ t\in[0,t^*-\varepsilon].\end{align}
This, together with (\ref{4.37}), shows that for each $t\in [0,
t^*-\varepsilon]$,
\begin{align}
\|\psi^\varepsilon(t)\|\leq&1-\|y^\varepsilon(t^*-\varepsilon)\|+\int_t^{t^*-
\varepsilon}\|f_y^{(2)}(y^\varepsilon(\tau))\|\ \|\psi^\varepsilon(\tau)\|d\tau\nonumber\\
\leq&1-\|y^\varepsilon(t^*-\varepsilon)\|+\int_t^{t^*-
\varepsilon}\frac{\sqrt{(1-\|y^\varepsilon(\tau)\|)^2+1}}{(1-\|y^\varepsilon(\tau))\|)^2}
\|\psi^\varepsilon(\tau)\|d\tau.\nonumber
\end{align}
Now we can apply  the Gronwall inequality to get that for each
$t\in [0, t^*-\varepsilon]$,
\begin{align}
\|\psi^\varepsilon(t)\|
\leq(1-
\|y^\varepsilon(t^*-\varepsilon)\|)\exp\Big\{
\int_t^{t^*-\varepsilon}\big(\frac{1}{(1-\|y^\varepsilon
(\tau)\|)^2}+\frac{1}{1-\|y^\varepsilon(\tau)\|}
\big)d\tau\Big\},\nonumber
\end{align}
from which and (\ref{4.38-1}),
it follows that for each
$t\in [0, t^*-\varepsilon]$,
\begin{align}
\label{4.40}\|\psi^\varepsilon(t)\|
\leq\frac{1-
\|y^\varepsilon(t^*-\varepsilon)\|}{\|y^\varepsilon(t^*-\varepsilon)\|}\exp\Big\{
\int_t^{t^*-\varepsilon}\big(\frac{1}{(1-\|y^\varepsilon
(\tau)\|)^2}+\frac{1}{1-\|y^\varepsilon(\tau)\|}
\big)d\tau\Big\}.
\end{align}

On the other hand, according to the system
satisfied by $y^\varepsilon(\cdot)$, it follows that for
every $t\in [0, t^*-\varepsilon]$,
\begin{align}
&\int_t^{t^*-\varepsilon}\frac{d\displaystyle\frac{1-
\|y^\varepsilon(\tau)\|}{\|y^\varepsilon(\tau)\|}}
{\displaystyle\frac{1-
\|y^\varepsilon(\tau)\|}{\|y^\varepsilon(\tau)\|}}=
\ln\frac{\displaystyle\frac{1- \|y^\varepsilon(
t^*-\varepsilon)\|}{\|y^\varepsilon(
t^*-\varepsilon)\|}}{\displaystyle\frac{1-
\|y^\varepsilon(t)\|}{\|y^\varepsilon(t)\|}}\nonumber\\=
&\int_t^{t^*-\varepsilon}
\Big\{-\frac{1}{(1-\|y^\varepsilon(\tau)\|)^2}-\frac{<y^\varepsilon
(\tau),B(\tau)u^\varepsilon(\tau)>
}{(1-\|y^\varepsilon(\tau)\|)(\|y^\varepsilon(\tau)\|)^2}\Big\}d\tau,\nonumber
\end{align}
namely, we have the equation as follows:
\begin{align}
&\displaystyle\frac{1-
\|y^\varepsilon(t^*-\varepsilon)\|}{\|y^\varepsilon(t^*-\varepsilon)\|}
\nonumber\\=&\displaystyle\frac{1-
\|y^\varepsilon(t)\|}{\|y^\varepsilon(t)\|}\cdot
\exp\Big\{\int_t^{t^*-\varepsilon}
\frac{-1}{(1-\|y^\varepsilon(\tau)\|)^2}d\tau\Big\}\nonumber\\&\cdot
\exp\Big\{\int_t^{t^*-\varepsilon} \frac{-<y^\varepsilon
(\tau),B(\tau)u^\varepsilon(\tau)>
}{(1-\|y^\varepsilon(\tau)\|)(\|y^\varepsilon(\tau)\|)^2}d\tau\Big\}, \ \
t\in [0, t^*-\varepsilon]. \nonumber\end{align} Then, making use of
the above equation and (\ref{4.40}), we obtain that for each
$t\in [0, t^*-\varepsilon]$,
\begin{align}
\label{4.41}\|\psi^\varepsilon(t)\| \leq\frac{1-
\|y^\varepsilon(t)\|}{\|y^\varepsilon(t)\|}\exp\Big\{
\int_t^{t^*-\varepsilon}\big(\frac{1}{1-\|y^\varepsilon(\tau)\|}
-\frac{<y^\varepsilon (\tau),B(\tau)u^\varepsilon(\tau)>
}{(1-\|y^\varepsilon(\tau)\|)(\|y^\varepsilon(\tau)\|)^2)}
\big)d\tau\Big\}.
\end{align} On the other hand, since $T_q(f^{(2)},y^0,u^\varepsilon)\geq t^*$, we can conclude from
Lemma 2.1.2 and Lemma 2.2.2 that for each $[0, t^*-
\varepsilon]$,
$$\|y^\varepsilon(t)\|>1-\frac{1}{2K_0+1},\
\frac{1}{1-\|y^\varepsilon(t)\|}\leq
C(T_q(f^{(2)},y^0,u^\varepsilon)-t)^{-{2}/{3}}\leq
C(t^*-t)^{-{2}/{3}}.\;\;\;
$$
This, together with (\ref{4.41}), implies that
\begin{align}\label{4.42}\|\psi^\varepsilon(t)\|\leq C_2(1-
\|y^\varepsilon(t)\|), \  \ t\in [0, t^*-\varepsilon],
\end{align}
where  $C_2$ is independent of $\varepsilon\in (0,t^*)$ and $t$.

Now, we can use the similar argument we used of the proof of Theorem 1.2 to find a
sequence $\{\varepsilon_{n,n}\}_{n=1}^\infty$ and a function $\psi(\cdot)\in C([0,t^*];\mathbb{R}^2)$
with $\psi(t^*)=0$ such that
\begin{align}\psi^{\varepsilon_{n,n}}(\cdot)\rightarrow
\psi(\cdot)\ \mbox{uniformly  on}\ [0, t^*-\delta]\ \
\mbox{as}\ n\rightarrow +\infty,\;\;\mbox{for each}\; \delta\in (0,
t^* )\nonumber
\end{align}
and
\begin{align}y^{\varepsilon_{n,n}}(\cdot)\rightarrow
y^*(\cdot)  \  \mbox{uniformly on}\ [0, t^*-\delta] \ \
\mbox{as}\ n\rightarrow +\infty,\;\;\mbox{for each}\; \delta\in (0,
t^* ).\nonumber
\end{align}

Next, we shall show the non-triviality of the function $\psi(\cdot)$ in
this case.

By the equations satisfied by $y^{\varepsilon_{n,n}}(\cdot)$ and
$\psi^{\varepsilon_{n,n}}(\cdot)$, respectively,  we see that
\begin{align}
&\int_t^{t^*-\varepsilon_{n,n}}\frac{d}{d\tau}<\frac{y^{\varepsilon_{n,n}}(\tau)}
{1-\|y^{\varepsilon_{n,n}}(\tau)\|},\psi^{\varepsilon_{n,n}}(\tau)>
d\tau\nonumber\\=
&\int_t^{t^*-\varepsilon_{n,n}}\Big\{<\frac{y^{\varepsilon_{n,n}}(\tau)}
{1-\|y^{\varepsilon_{n,n}}(\tau)\|},
-f_y^{(2)}(y^{\varepsilon_{n,n}}(\tau))\psi^{\varepsilon_{n,n}}(\tau)>\nonumber\\&
+<[f_y^{(2)}(y^{\varepsilon_{n,n}}(\tau))]^T\Big(\frac{y^{\varepsilon_{n,n}}(\tau)}
{1-\|y^{\varepsilon_{n,n}}(\tau)\|}+B(\tau)u^{\varepsilon_{n,n}}(\tau)\Big),
\psi^{\varepsilon_{n,n}}(\tau)>\Big\}d\tau\nonumber\\=&\int_t^{t^*-\varepsilon_{n,n}}
<B(\tau)u^{\varepsilon_{n,n}}(\tau),
f_y^{(2)}(y^{\varepsilon_{n,n}}(\tau))\psi^{\varepsilon_{n,n}}(\tau)>d\tau,
\ t\in [0, t^*\mbox{$-$}{\varepsilon_{n,n}}],\nonumber
\end{align}
from which and (\ref{4.38}), it follows that for $\ t\in [0,
t^*-\varepsilon_{n,n}]$,
\begin{align}
&\|y^{\varepsilon_{n,n}}(t^{*}-\varepsilon_{n,n})\|-<\frac{y^{\varepsilon_{n,n}}(t)}
{1-\|y^{\varepsilon_{n,n}}(t)\|},\psi^{\varepsilon_{n,n}}(t)>\nonumber\\
=&\int_t^{t^*-\varepsilon_{n,n}}
<B(\tau)u^{\varepsilon_{n,n}}(\tau),
f_y^{(2)}(y^{\varepsilon_{n,n}}(\tau))\psi^{\varepsilon_{n,n}}(\tau)>d\tau.\nonumber
\end{align}
Then, let $n\rightarrow +\infty$ in the above equation, it holds that for each $t\in [0,t^*)$,
\begin{align}
&1-<\frac{y^*(t)}
{1-\|y^*(t)\|},\psi(t)>\nonumber\\
=&\int_t^{t^*} <B(\tau)u^*(\tau),
f_y^{(2)}(y^*(\tau))\psi(\tau)>d\tau.\nonumber
\end{align}
This shows that the function $\psi(\cdot)$ is not trivial.

Finally, we can use  the similar argument we used of the proof of Theorem 1.2 to
prove that $\psi(\cdot)$ satisfies (\ref{1.5}) and (\ref{1.6}).
The inequality (\ref{1.6-1}) can be proved by (\ref{2.22-1}) in Lemma
2.2.2.

This completes the proof of Theorem 1.3.\ \ \ \#\\

{\it The author gratefully acknowledges Professor Gengsheng Wang for his useful
suggestions and help.}

\end{document}